\documentclass[reqno,11pt]{amsart}

\usepackage{amsfonts,amsmath,amsthm,amssymb}
\usepackage{shortvrb,graphicx,psfrag}
\usepackage{wrapfig}
\usepackage[colorlinks,pdfpagelabels,pdfstartview = FitH,bookmarksopen
= true,bookmarksnumbered = true,linkcolor = blue,plainpages =
false,hypertexnames = false,citecolor = red,pagebackref=false]{hyperref}
\usepackage{relsize}
\usepackage{xcolor}
\usepackage{appendix}
\usepackage{comment}
\usepackage{ulem}
\usepackage[makeroom]{cancel}
\usepackage[margin=1.4in]{geometry}
\usepackage{stmaryrd}
\usepackage{float}
\SetSymbolFont{stmry}{bold}{U}{stmry}{m}{n}

\allowdisplaybreaks

\let\TeXchi\chi
\newbox\chibox
\setbox0 \hbox{\mathsurround0pt $\TeXchi$}
\setbox\chibox \hbox{\raise\dp0 \box 0 }
\def\chi{\copy\chibox}

\newcommand{\N}{\mathbb N}
\renewcommand{\d}{\mathrm{d}}
\newcommand{\dx}{\mathrm{d}x}

\newcommand{\dt}{\mathrm{d}t}
\newcommand{\ds}{\mathrm{d}s}
\renewcommand{\rho}{\varrho}

\newcommand{\datap}{\{N,p,C_o,C_1\}}
\newcommand{\bry}{B_{\rho}(y)}
\newcommand{\kry}{K_{\rho}(y)}

\newcommand{\noi}{\noindent}
\newcommand{\dsty}{\displaystyle}
\newcommand{\txty}{\textstyle}
\newcommand{\pl}{\partial}

\newcommand{\al}{\alpha}

\newcommand{\be}{\beta}

\newcommand{\gm}{\gamma}
\newcommand{\dl}{\delta}
\newcommand{\Dl}{\Delta}
\newcommand{\lm}{\lambda}

\newcommand{\varep}{\varepsilon}
\newcommand{\eps}{\epsilon}

\newcommand{\sig}{\sigma}

\newcommand{\om}{\omega}
\newcommand{\Om}{\Omega}

\newcommand{\z}{\zeta}

\newcommand{\loc}{\operatorname{loc}}

\newcommand{\dist}{\operatorname{dist}}
\newcommand{\dvg}{\operatorname{div}}
\newcommand{\essup}{\operatornamewithlimits{ess\,sup}}
\newcommand{\essinf}{\operatornamewithlimits{ess\,inf}}
\newcommand{\osc}{\operatornamewithlimits{osc}}
\newcommand{\bl}[1]{\mathbf{#1}}
\newcommand{\rr}{\mathbb{R}}
\newcommand{\rn}{\rr^N}
\newcommand{\nn}{\mathbb{N}}
\newcommand{\df}[1]{\buildrel\mbox{\small def}\over{#1}}

\def\Xint#1{\mathchoice
    {\XXint\displaystyle\textstyle{#1}}
    {\XXint\textstyle\scriptstyle{#1}}
    {\XXint\scriptstyle\scriptscriptstyle{#1}}
    {\XXint\scriptscriptstyle\scriptscriptstyle{#1}}
    \!\int}
\def\XXint#1#2#3{\setbox0=\hbox{$#1{#2#3}{\int}$}
    \vcenter{\hbox{$#2#3$}}\kern-0.5\wd0}
\def\bint{\Xint-}
\def\dashint{\Xint{\raise4pt\hbox to7pt{\hrulefill}}}
\def\dashiint{\bint\kern-0.15cm\bint}

\newtheorem{proposition}{Proposition}[section]
\newtheorem{theorem}{Theorem}[section]

\newtheorem{lemma}{Lemma}[section]

\newtheorem{remark}{Remark}[section]
\newtheorem{definition}{Definition}[section]

\numberwithin{equation}{section}
\numberwithin{theorem}{section}
\numberwithin{proposition}{section}
\numberwithin{lemma}{section}
\numberwithin{remark}{section}
\title[Improved moduli of continuity for degenerate phase transitions]{Improved moduli of continuity\\ for degenerate phase transitions}

\author[U. Gianazza]{Ugo Gianazza}
\address{Dipartimento di Matematica ``F. Casorati", 
Universit\`a di Pavia, via Ferrata 5, 27100 Pavia, Italy}{}
\email{ugogia04@unipv.it}

\author[N. Liao]{Naian Liao}
\address{Fachbereich Mathematik, Universit\"at Salzburg,
Hellbrunner Str. 34, 5020 Salzburg, Austria}{}
\email{naian.liao@plus.ac.at}

\author[J.M.~Urbano]{Jos\'{e} Miguel Urbano}
\address{Applied Mathematics and Computational Sciences Program (AMCS), Compu\-ter, Electrical and Mathematical Sciences and Engineering Division (CEMSE), King Abdullah University of Science and Technology (KAUST), Thuwal, 23955 -6900, Kingdom of Saudi Arabia and CMUC, Department of Mathematics, University of Coimbra, 3000-143 Coimbra, Portugal}{} 
\email{miguel.urbano@kaust.edu.sa}

\begin{document}



\begin{abstract}
We substantially improve in two scenarios the current state-of-the-art modulus of continuity for weak solutions to the $N$-dimensional, two-phase Stefan problem featuring a $p-$degenerate diffusion: for $p=N\geq 3$, we sharpen it to
$$
\boldsymbol{\omega}(r) \approx \exp (-c| \ln r|^{\frac1N});
$$
for $p>\max\{2,N\}$, we derive an unexpected H\"older modulus. 
\end{abstract}  

\date{\today}

\maketitle

\section{Introduction}

The classical two-phase Stefan problem is an archetypal free boundary problem that models a phase transition at constant temperature. It consists of solving the heat equation (or nonlinear variants of it) in the solid and liquid phases, coupled with the so-called Stefan condition at the \textit{a priori} unknown interface separating them. This condition corresponds to an energy balance, prescribing the proportionality between the jump of the heat flux across the free boundary and its local velocity.

In this paper, we are interested in a phase transition problem involving a degenerate diffusion of $p-$Laplacian-type, for $p>2$. In its weak form, for which any explicit reference to the free boundary is absent, such a problem can be formulated as
\begin{equation}\label{Eq:1:1}
\pl_t\be(u)-\dvg\bl{A}(x,t,u, Du) \ni 0\quad \text{ weakly in }\> E_T,
\end{equation}
for a function $u:E_T \to \rr$ representing the temperature.
Here, $E_T\df=E\times(0,T]$, for some open set $E\subset \rn$, $N\in\nn$ and $T>0$,
whereas $\be(\cdot)$ is the maximal monotone graph defined by
\begin{equation}\label{Eq:beta}
\be(s)=\left\{
\begin{array}{cl}
s\quad& \text{if}\> s>0,\\[5pt] 
[-\nu,0]\quad&\text{if}\> s=0,\\[5pt]
s-\nu\quad&\text{if}\> s<0,
\end{array}
\right.
\end{equation} 
for a positive constant $\nu$, the latent heat of the phase transition, representing the aforementioned proportionality ratio.   
The function $\bl{A}(x,t,z,\xi)\colon E_T\times \rr \times \rr^{N}\to\rr^N$ is a {\it Carath\'eo\-dory} function, satisfying the structure conditions
\begin{equation}  \label{Eq:1:2}
\left\{
\begin{array}{l}
\bl{A}(x,t,z,\xi)\cdot \xi\ge C_o|\xi|^p \\[5pt]
|\bl{A}(x,t,z,\xi)|\le C_1|\xi|^{p-1} 
\end{array}
\right .\quad \text{ a.e.}\> (x,t)\in E_T,\, \forall\,z\in\rr,\,\forall\,\xi\in\rr^N,
\end{equation}
where $0<C_o\le C_1$ are given positive constants. 

Contrary to what was a common belief at the inception of the modern mathematical theory of the problem, the temperature $u$ turns out to be a continuous function. This was proven independently by Caffarelli-Evans \cite{Caff-Ev}, DiBenedetto \cite{DB-82}, Sacks \cite{Sacks}, and Ziemer \cite{Ziemer} in the early 1980s, much to the surprise of the Russian school that developed the weak formulation and the corresponding existence theory in the 1960s (see \cite{Olienik}, and also \cite{Visintin}). The result would eventually be extended in \cite{Urb} to the case $p>2$ treated here. 

The quantification of these qualitative results, \textit{i.e.}, the search for a modulus of continuity for the temperature $u$, has a history of its own. Besides its intrinsic interest, the knowledge of a quantitative modulus of continuity is important for applications. Indeed, if we define the (weak) free boundary as $\partial[u>0]$, there is a direct connection between the regularity of this interface and the modulus of continuity of $u$, as discussed, for example, in \cite[Section~6]{Caff-Fried}.

A modulus of the type 
$$
\boldsymbol{\omega}(r) \approx \left[ \ln \left| \ln \left( \frac{r}{R}  \right) \right| \right]^{-\sigma}, \qquad  r\in\left(0,R\right),
$$
for some $\sigma >0$, which follows from the reasoning in \cite{DB-82, Caff-Ev}, would eventually appear explicitly in \cite[Remark 3.1]{DiB-Fried} and \cite{DB-86}. It remained the state-of-the-art for 30 years until the next decisive result in \cite{BKU-14} improved the modulus to 
\begin{equation*}  
\boldsymbol{\omega}(r) \approx \left| \ln \left( \frac{r}{R}  \right) \right|^{-\sigma (N,p)}, \qquad  r\in\left(0,R\right),
\end{equation*}
including the degenerate case $p>2$, discarding an iteration of the logarithm and determining the precise value of the exponent $\sigma$ in terms of the data of the problem: in the case $2<p<N$, for example, $\sigma = \frac{p}{N+p}$. The conjecture stated in \cite{BKU-14} that the sharp, optimal modulus of continuity for the two-phase Stefan problem had been obtained turned out to be speculative, and we disprove it here, obtaining a better modulus in two scenarios: for $p = N$ in dimensions $N=3,4, \ldots$, the modulus is upgraded to 
\[
\boldsymbol{\omega}(r) \approx \exp\left(-c\left|\ln\left(\frac{r}{R}\right)\right|^{\frac1N}\right), \qquad  r\in\left(0,R\right);
\]
for dimensions $N=1,2$, and for $p > N$ in dimensions $N=3,4, \ldots$, the modulus is improved to a H\"older modulus of type
\[
\boldsymbol{\omega}(r) \approx \left(\frac rR\right)^{\gamma}, \qquad  r\in\left(0,R\right).
\]
These findings constitute the object of our main result.

\begin{theorem}\label{Thm:1}
Let $u$ be a bounded, local weak solution to the Stefan problem \eqref{Eq:1:1}--\eqref{Eq:1:2} and let  $$
\mathbb M=\max\Big\{1,\osc_{E_T}u\Big\}.
$$ Consider a compact subset $\mathcal K$ of $E_T$ and denote by $R$ the $p-$parabolic distance of $\mathcal K$ to the parabolic boundary of $E_T$, cf.~\eqref{Eq:par-dist}. 

\medskip

\noi\textbf{(I)}~When $p=N\ge3$, there exist constants $c,\,\sig\in(0,1)$, depending on $\{\nu, N, C_o, C_1\}$, such that for every $(x_i,t_i)\in\mathcal{K}$, $i=0,1$, satisfying 
\[
|x_o-x_1|+\mathbb M^{\frac{N-2}{N}}|t_o-t_1|^{\frac1N}\le\sig R,
\] we have
\[
|u(x_o,t_o)-u(x_1,t_1)|\le \mathbb M \exp\left(-c\left|\ln\left(\frac{|x_o-x_1|+\mathbb M^{\frac{N-2}{N}}|t_o-t_1|^{\frac1N}}{\sig R} \right)\right|^{\frac1N}\right).
\]

\noi\textbf{(II)}~When $p>\max\{2,N\}$, there exist constants $\gm\in(0,1)$ and $C>1$, depending on $\{N, p, C_o, C_1\}$ but independent of $\nu$, such that for every $(x_i,t_i)\in\mathcal{K}$, $i=0,1$, satisfying 
\[
|x_o-x_1|+\mathbb M^{\frac{p-2}{p}}|t_o-t_1|^{\frac1p}\le R,
\] we have
\[
|u(x_o,t_o)-u(x_1,t_1)|\le C\,\mathbb M  \left(\frac{|x_o-x_1|+\mathbb M^{\frac{p-2}{p}}|t_o-t_1|^{\frac1p}}{R} \right)^\gm .
\] 
\end{theorem}

\vskip.2truecm

For the remaining case, namely $2<p<N$, $N\ge3$, our approach yields a logarithmic type modulus of continuity like the one obtained in \cite{BKU-14}. Therefore, the significance of our contribution lies in bringing previously hidden structural properties of solutions to light and providing a unified perspective on the issue in all scenarios. Table~\ref{Table:1} summarises the current state-of-the-art in terms of the modulus of continuity (MoC) for the two-phase Stefan problem involving a degenerate diffusion of $p-$Laplacian-type, for $p>2$.

\begin{table}[ht]
\centering
\begin{tabular}{p{4.3cm}p{4.3cm}p{1cm} p{1cm}|p{1cm}|p{1cm}|p{1cm}|p{1.4cm}|}
\hline
\qquad $p$ \& $N$ & \qquad MoC \\ \hline
$p>\max\{2,N\}$ & $\boldsymbol\om^{(1)}(r)\approx r^\gamma$ \\ 
$p=N\ge3$ & $\boldsymbol\om^{(2)}(r)\approx \exp(-c|\ln r|^{\frac1N})$\\
$2<p<N$, $N\ge3$ & $\boldsymbol\om^{(3)}(r)\approx|\ln r|^{-\sig}$ \\
\hline
\end{tabular}
\vspace{.2cm}
\caption{Moduli of Continuity}\label{Table:1}
\end{table}

Among these three types of moduli, the H\"older modulus of continuity $\boldsymbol\om^{(1)}(\cdot)$ is perhaps the most surprising one. Loosely speaking, it shows the $p-$Laplacian-type diffusion dictates the local behaviour of solutions as though the singularity of $\be(\cdot)$ plays no role. The H\"older exponent $\gamma$ determined by our method is implicit (see \cite{Liao-p>N} for a slightly different approach), and it would be interesting to search for an explicit $\gamma$, even in the archetypal case of the $p-$Laplacian. On the other hand, the modulus $\boldsymbol\om^{(2)}(\cdot)$ represents a sort of borderline case; it was observed by Caffarelli and Friedman in their study \cite{Caff-Fried} of the so-called ``one-phase" Stefan problem and, recently, established in \cite{Liao-2D}, in the case $p=N=2$, under a general framework similar to that of this paper.

\subsection{New techniques}
The singularity of $\beta$ occurs at $u=0$; henceforth, it is natural to expect that the set $[u=0]$ and its measure play a fundamental role in the derivation of the modulus of continuity. As is frequently the case when studying the regularity of solutions to the two-phase Stefan problem, the need to deal with such a measure leads to considering an alternative. However, a significant difference appears here: instead of taking into account the measure of $[u=0]$ inside a proper cylinder, we deal with measures inside time slices, that is, sets of the form $[u(\cdot,t)=0]\cap K_\rho$, for $t$ in a suitable interval $I$, where $K_\rho$ denotes the cube of wedge $2\rho$. Therefore, the alternative is stated as follows: either there exists a $\bar t\in{I}$ such that $\left|[u(\cdot,\bar t)=0] \right|$ is small in $K_\rho$, or for every $t\in{I}$, the set where $u(\cdot,t)$ is close to zero occupies a sizable portion of the cube. 

The former possibility is the favourable one, and it is dealt with in Sections~\ref{SS:4:1} and \ref{SS:4:3}; roughly speaking, the smallness of the singularity set yields that the solution behaves as if it were the solution of a $p-$Laplacian-type equation, and its continuity can be obtained in a fairly standard way.

The latter occurrence is more complicated, and it is tackled in Sections~\ref{SS:4:2} and \ref{SS:4:4} thanks to Proposition~\ref{Prop:3:1}, which is a new {\it expansion of positivity} for non-negative super-solutions to parabolic $p-$Laplacian-type equations as \eqref{Eq:3:1}, satisfying the structure conditions \eqref{Eq:1:2}. It is precisely this result that allows us to discriminate the behaviour of solutions in terms of the dimension $N$. We believe it might be of independent interest. 

To better understand the novelty, let us first consider the by-now classical expansion of positivity for such a super-solution $v$ to parabolic $p-$Laplacian-type equations when $p>2$ (see \cite[Chapter~5, Proposition~7.1]{DBGV-mono}): if, for some $k>0$ and $\alpha\in(0,1)$,
\[
|[v(\cdot,0)\ge k]\cap K_1|\ge\alpha|K_1|,
\]
there exist constants $\xi \in(0,1)$, and $b,\,d>1$, depending on $\{N,p,C_o,C_1\}$, such that
\[
v(\cdot,t)\ge\xi \alpha^d k\qquad\text{ in }\,\, K_2,
\]
for all times
\[
\frac{b}{(\xi\alpha^d k)^{p-2}} \le t\le \frac{2b}{(\xi\alpha^d k)^{p-2}}.
\]
The crucial issue is the power-like dependence of the shrinking parameter $\xi \alpha^d$ on the measure-theoretical quantity $\alpha$. In contrast, Proposition~\ref{Prop:3:1} yields that, given a super-solution $v$, if 
\[
|[v(\cdot,t)\ge k]\cap K_1|\ge\alpha|K_1|,
\]
for all $t\in(0,(\bar\delta\,{\rm cap}_p(K_{\alpha^2}, K_{3})\,k)^{2-p}]$, for a proper $\bar\delta$ which depends only on the data $\{N,p,C_o,C_1\}$, there exist constants $\xi\in(0,1)$ and $b>1$, with the same dependence as $\bar\delta$, such that
\[
v(\cdot,t)\ge\xi\, {\rm cap}_p(K_{\alpha^2}, K_{3})\, k\qquad\text{ in }\,\, K_2,
\]
for all times
\[
\frac{b}{(\xi\,{\rm cap}_p(K_{\alpha^2}, K_{3})\,k)^{p-2}}\le t\le \frac{2b}{(\xi\,{\rm cap}_p(K_{\alpha^2}, K_{3})\,k)^{p-2}}.
\]
The different sets where the initial measure-theoretical lower bound is achieved are just technical differences and do not play any substantial role; the main point is that the shrinking parameter changes from $\xi \alpha^d$ to $\xi\, {\rm cap}_p(K_{\alpha^2}, K_{3})$ and, relying on the properties of the $p-$capacity, we can dramatically improve the dependence with respect to $\alpha$ when $p\ge N$, and consequently obtain the sharpened moduli $\boldsymbol\om^{(1)}$ and $\boldsymbol\om^{(2)}$ of Table~\ref{Table:1}.
The improvement of the shrinking parameter is based on the deep measure-theoretical Lemma~\ref{Lm:0:3:1}, and a proper adaptation of the capacity estimates originally employed in \cite{GL-deg}. A result similar to Proposition~\ref{Prop:3:1} was first used in \cite{Liao-2D} for non-degenerate parabolic equations; extending such a statement to the degenerate context $p>2$ is far from trivial.
Indeed, it requires careful use of the sophisticated intrinsic scaling method, typical of parabolic $p-$Laplacian-type equations, for which we refer to \cite{dibe-sv, DBGV-mono, Urb-08}.

Finally, it is also worth mentioning that our approach forgoes the use of logarithmic estimates, unlike most previous works on the topic. The mechanism we employ to spread information in time is based on Lemma \ref{Lm:DG:initial:1}, which follows from energy estimates alone.

\subsection{Organization of the paper}
We organize the paper as follows. In Section~\ref{S:prelim}, we introduce the basic notation and definitions and collect some auxiliary results. In Section~\ref{S:3}, we present the new expansion of positivity and its proof, whereas the proof of the main theorem is given in Section~\ref{S:proof-thm}.

\section{Preliminary results}\label{S:prelim}

This section gathers some useful pieces of notation, the main definitions, and several auxiliary results that will be instrumental in the main proofs. They comprise a measure-theoretical lemma, energy estimates and three De Giorgi-type lemmata; for the third, a detailed proof is included since it generalizes to any $p>1$ the corresponding result for $p=2$. 

\subsection{Notation and definitions}
For $\rho>0$ and $y\in\rn$, denote by $\kry$ the cube of edge $2\rho$, centred at $y$ and with faces 
parallel to the coordinate planes, and by $\bry$ the ball of radius $\rho$ centred at $y$. If $y$ is the origin, 
let $K_\rho(0)=K_\rho$, and $B_\rho(0)=B_\rho$. Most of the statements are given in terms of cubes, but we briefly refer to balls in Appendix~\ref{A:1} when considering simple, explicit examples of $p-$capacity.

The parabolic boundary of the cylinder $E_T$ is given by the union of its lateral and initial boundary
\[
\partial_{\textrm{par}}E_T\df=(\overline E\times\{0\})\cup(\partial E\times(0,T]).
\]
For any compact subset $\mathcal K$ of $E_T$, the $p-$parabolic distance of $\mathcal K$ to the parabolic boundary $\partial_{\textrm{par}}E_T$ is defined by
\begin{equation}\label{Eq:par-dist}
p-\dist_{\textrm{par}}({\mathcal K},\partial_{\textrm{par}}E_T)\df=\inf_{\genfrac{}{}{0pt}{}{(x,t)\in{\mathcal K}}{(y,s)\in\partial_{\textrm{par}}E_T}}\left(\max\left\{\max_{1\le i\le N}|x_i-y_i|,\mathbb M^{\frac{p-2}p}|t-s|^{\frac1p}\right\}\right).
\end{equation}
As is usual with structure conditions like those in \eqref{Eq:1:2}, the definition of $p-$parabolic distance depends on the solution itself. It is not hard to see that the previous definition is equivalent to the one usually employed with the parabolic $p-$Laplacian, namely
\[
p-\dist_{\textrm{par}}({\mathcal K},\partial_{\textrm{par}}E_T) = \inf_{\genfrac{}{}{0pt}{}{(x,t)\in{\mathcal K}}{(y,s)\in\partial_{\textrm{par}}E_T}}\left(|x-y|+\mathbb M^{\frac{p-2}p}|t-s|^{\frac1p}\right).
\]
Using \eqref{Eq:par-dist} simplifies some of the computations at the end of this work (cf. Subsection \ref{4.6}).

For $R,\,S,\,\rho,\,\theta>0$, the backward and forward cylinders with vertex $(x_o,t_o)$ in $\rr^{N+1}$ are defined by
\begin{equation*}
\left\{
\begin{array}{ll}
(x_o,t_o)+Q_{R,S}\df=K_R(x_o)\times (t_o-S,t_o],\\[5pt]
(x_o,t_o)+Q_{\rho}(\theta)\df=K_{\rho}(x_o)\times(t_o-\theta\rho^p,t_o],\\[5pt]
(x_o,t_o)+Q^+_{\rho}(\theta)\df=K_{\rho}(x_o)\times(t_o,t_o+\theta\rho^p].
\end{array}
\right.
\end{equation*}
Functions truncated by $k\in\rr$ are defined as
\[
(u-k)_+\df=\max\{u-k, 0\},\qquad (u-k)_-\df=\max\{k-u, 0\}.
\] 
For a Lebesgue measurable set $A\subset\rn$, we denote its measure by $|A|$. In various estimates, we use $\gm$ to represent a generic positive constant.

\medskip

We now introduce the precise notions of local weak sub-solution, local weak super-solution, and local weak solution, which we will deal with throughout the paper. The relationship between classical and weak formulations of the two-phase Stefan problem is an interesting issue in itself, for which we refer to \cite{Visintin}.

\begin{definition}
A function
\begin{equation*}
	u\in L_{\loc}^{\infty}\left(0,T;L^2_{\loc}(E)\right)\cap L^p_{\loc}\left(0,T; W^{1,p}_{\loc}(E)\right)
\end{equation*}
is a local weak sub(super)-solution to the Stefan problem~\eqref{Eq:1:1} with the structure
conditions \eqref{Eq:1:2}, if, for every compact set $\mathcal K \subset E$ and every sub-interval
$[t_1,t_2]\subset (0,T]$, there is a selection $v\subset\be(u)$, \textit{i.e.},
\[
\left\{\left(z,v(z)\right): z\in E_T\right\}\subset \left\{\left(z,\be[u(z)]\right): z\in E_T\right\},
\]
 such that
\begin{equation*}
	\int_{\mathcal K} v\z \,\dx\bigg|_{t_1}^{t_2}
	+
	\iint_{\mathcal K \times(t_1,t_2)} \left[-v\pl_t\z+\bl{A}(x,t,u,Du)\cdot D\z\right]\dx\dt
	\le(\ge)0,
\end{equation*}
for all non-negative testing functions
\begin{equation*}
\z\in W^{1,2}_{\loc}\left(0,T;L^2(\mathcal K)\right)\cap L^p_{\loc}\left(0,T;W_o^{1,p}(\mathcal K)
\right).
\end{equation*}
A function that is both a local weak sub-solution and a local weak super-solution is termed a local weak solution.
\end{definition}

The above notion of solution, though natural from the viewpoint of existence theory, renders an issue when one attempts to use the solution as a testing function since the time derivative does not exist in the Sobolev sense. To confront it,
we will consider the regularized version of the Stefan problem \eqref{Eq:1:1}:
\begin{equation}\label{Eq:reg}
\begin{aligned}
&\pl_t \be_\eps(u) -\dvg\bl{A}(x,t,u, Du) = 0,\quad \text{ weakly in }\> E_T.
\end{aligned}
\end{equation}
Here, we define the mollification of $\be$ by
\begin{equation*}
\be_\eps(s)\df=s+\nu H_{\eps}(s),
\end{equation*}
for a parameter $\eps\in(0,1)$ and
\begin{equation*}
H_{\eps}(s)\df=\left\{
\begin{array}{cc}
0,\quad& s>0,\\[5pt] 
\frac1{\eps}s,\quad& -\eps\le s\le 0,\\[5pt]
-1,\quad& s<-\eps.
\end{array}
\right.
\end{equation*} 

The notion of solution for \eqref{Eq:reg} can be found in~\cite[Chapter~II]{dibe-sv}. 
We will proceed with the assumption that local solutions to the Stefan problem \eqref{Eq:1:1}
can be approximated by a sequence of solutions to \eqref{Eq:reg} locally uniformly. This is a standard assumption, cf.~\cite{DBV-95, BKU-14, GL-Stefan-1}.

It is well-known that, for any fixed $\eps$, solutions to \eqref{Eq:reg} are locally H\"older continuous; see~\cite[Chapter~III, IV]{dibe-sv}. However, such H\"older estimates will not be retained as $\eps\to0$. We aim to identify estimates on the modulus of continuity that are stable as $\eps\to0$. In this way, the Stefan problem \eqref{Eq:1:1}, as a limiting case of the regularized problem~\eqref{Eq:reg}, enjoys the same kind of modulus of continuity. This kind of regularizing scheme also appears in the study of gradient regularity for congested traffic dynamics \cite{BDGP-1, BDGP-2}.

\subsection{A measure-theoretical lemma}
The following result exhibits the clustering of positivity for Sobolev functions in a measure-theoretical sense.

\begin{lemma}\label{Lm:0:3:1}
Let $u\in W^{1,1}(K_\rho)$ satisfy 
\begin{equation*}
\bint_{K_\rho}|Du|\,\dx\le\gm\frac{k}{\rho}\quad\text{ and }\quad 
|[u>k]\cap K_{\rho}|\ge\al|K_\rho|,
\end{equation*}
for some positive $k,\,\gm$, and $\al\in(0,1)$. Then, for 
every $\dl,\,\lm\in(0,1)$, there exist $\varep=\varep(\al,\dl,\gm,\lm,N)\in(0,1)$ and $K_{\varep\rho}(y)\subset K_\rho$, such that 
\begin{equation*}
|[u>(1-\lm)k]\cap K_{\varep\rho}(y)|>(1-\dl)|K_{\varep\rho}|.
\end{equation*}
\end{lemma}

\begin{remark}\label{Rmk:0:3:1} \upshape
Following the various steps of the proof (for which we refer to \cite[Chapter 2, Lemma~3.1]{DBGV-mono}), the dependence 
of the reducing parameter $\varep$ on the measure-theoretical 
parameter $\al$, and on the constant $\gm$ appearing 
in the assumptions of the lemma can be traced to be 
of the form
\begin{equation}\label{Eq:0:3:7}
\varep=\bar{c}\,{\al^2}
\end{equation}
for a constant $\dsty\bar{c}=\frac{\lm \dl}{\gm C}$, where $C>1$ depends only on $N$, and is  
independent of $\al$.
\end{remark} 

\subsection{Energy estimates}

Consider a reference cylinder $(y,\tau)+Q_{R,S}\subset E_T$. For simplicity, we omit the reference point in the sequel and denote the cylinder as $Q_{R, S}$.
The following result has been proven in \cite[Proposition~2.1]{GL-Stefan-1}. 

\begin{proposition}\label{Prop:2:1}
	Let $p>1$ and consider a local weak sub(super)-solution $u$ to the Stefan problem~\eqref{Eq:reg} with structure conditions~\eqref{Eq:1:2} in $E_T$.
	There exists a positive constant $\gm=\gm (p,C_o,C_1)$, such that
 	for all cylinders $Q_{R,S}\subset E_T$,
 	every $k\in\rr$, and every non-negative, piecewise smooth cutoff function
 	$\z$ vanishing on $\pl K_{R}(y)\times(\tau-S,\tau]$,  there holds
\begin{align}\label{Eq:en-est-gen}
	\essup_{\tau-S<t<\tau}&\frac12\int_{K_R(y)\times\{t\}}	
	\z^p (u-k)_\pm^2\,\dx\nonumber\\
	&\quad + \nu\essup_{\tau-S<t<\tau}\int_{K_R(y)\times\{t\}}	
	\left(\int_0^{(u-k)_\pm}H'_\eps(k\pm s) s\,\ds\right)\z^p\,\dx\nonumber\\
	&\quad+
	\iint_{Q_{R,S}}\z^p|D(u-k)_\pm|^p\,\dx\dt\nonumber\\
	&\le
	\gm\iint_{Q_{R,S}}
		\left[
		(u-k)^{p}_\pm|D\z|^p + (u-k)_\pm^2|\partial_t\z^p|
		\right]
		\,\dx\dt\\
	&\quad+\nu \iint_{Q_{R,S}}\left(\int_0^{(u-k)_\pm}H'_\eps(k\pm s)s\,\d s\right)|\partial_t\z^p|\,\dx\dt\nonumber\\
	&\quad
	+\frac12\int_{K_R(y)\times \{\tau-S\}} \z^p (u-k)_\pm^2\,\dx\nonumber\\
	&\quad+\nu \int_{K_R(y)\times\{\tau-S\}}	
	\left(\int_0^{(u-k)_\pm}H'_\eps(k\pm s) s\,\ds\right)\z^p\,\dx\nonumber.
\end{align}
\end{proposition}
Here, we adopted the convention of interpreting $\pm$ as $+$ for sub-solutions and as $-$ for super-solutions. This convention also applies to the statements to follow.

\subsection{De Giorgi-type lemmata}
For $Q_{R,S} \subset E_T$, introduce the numbers $\mu^{\pm}$ and $\om$ satisfying
\begin{equation*}
	\mu^+\ge\essup_{Q_{R,S}} u,
	\quad 
	\mu^-\le\essinf_{Q_{R,S}} u,
	\quad
	\om\ge\mu^+ - \mu^-.
\end{equation*}
We discuss three lemmata resulting from the energy estimate~\eqref{Eq:en-est-gen} and De Giorgi-type iterations. It is worth mentioning that they all hold for $p>1$.

The first De Giorgi-type lemma is the same as \cite[Lemma~2.1]{GL-Stefan-1}; see~\cite[Lemma~2.1]{Liao-Stefan} for a proof.
\begin{lemma}\label{Lm:DG:1}
Let $u$ be a local weak sub(super)-solution to the Stefan problem~\eqref{Eq:reg} with structure conditions~\eqref{Eq:1:2} in $E_T$. Let $\delta, \,\xi\in(0,1)$, and set $\theta=\dl(\xi\om)^{2-p}$.
There exists a constant $c_o\in(0,1)$, depending only on the data $\{\nu, N, p, C_o, C_1\}$, such that if $\xi\om\le 1$ and if
\[
|[\pm(\mu^\pm -u)\le\xi\om]\cap [(x_o,t_o)+Q_{\rho}(\theta)] |\le c_o \delta^{\frac{N}{p}} (\xi\om)^{\frac{N+p}p}|Q_{\rho}(\theta)|,
\]
then
\[
\pm(\mu^\pm -u)\ge\tfrac12\xi\om, \quad\text{ a.e. in } (x_o,t_o)+Q_{\frac12\rho}(\theta),
\]
provided the cylinder $(x_o,t_o)+Q_{\rho}(\theta)$ is included in $Q_{R,S}$. 
\end{lemma}

\begin{remark}\label{Rmk:c_o}\upshape
Tracing the various constants in the proof gives that $c_o=\widetilde c_o\nu_\ast^{-(N+p)/p}$, where $\nu_\ast = \max\{1,\nu\}$ and $\widetilde c_o\in(0,1)$ depends only on $\{N,p, C_o, C_1\}$.
\end{remark}

The second De Giorgi-type lemma can be retrieved from \cite[Lemma~2.2]{GL-Stefan-1} or \cite[Lemma~2.1]{Liao-p>N}.
\begin{lemma}\label{Lm:DG:initial:1}
Let $u$ be a local weak sub(super)-solution to the Stefan problem~\eqref{Eq:reg} with structure conditions~\eqref{Eq:1:2} in $E_T$. 
Assume that, for some $\xi\in(0,1)$, there holds
\[
\pm(\mu^\pm -u(\cdot, t_1))\ge \xi\om,\quad\text{ a.e. in } K_{\rho}(x_o).
\]
There exists a constant $\gm_o\in(0,1)$, depending only on the data $\{ N, p, C_o, C_1\}$ and independent of $\nu$, such that, for any $\theta>0$, if 
\[
|[\pm(\mu^\pm -u)\le\xi\om]\cap [(x_o, t_1)+Q^+_{\rho}(\theta)]|\le\frac{ \gm_o(\xi\om)^{2-p}}{\theta}|Q^+_{\rho}(\theta)|,
\]
then
\[
\pm(\mu^\pm -u )\ge\tfrac12\xi\om,\quad\text{ a.e. in }K_{\frac12\rho}(x_o)\times(t_1,t_1+\theta\rho^p],
\]
provided the cylinder $(x_o, t_1)+Q^+_{\rho}(\theta)$ is included in $Q_{R,S}$. 
\end{lemma}

The third De Giorgi-type lemma generalizes to any $p>1$ the corresponding result for $p=2$ obtained in \cite[Lemma~2.2]{Liao-2D}. 
\begin{lemma}\label{Lm:DG-new}
Let $u$ be a local weak sub(super)-solution to the Stefan problem~\eqref{Eq:reg} with structure conditions~\eqref{Eq:1:2} in $E_T$.
For $\xi\in(0,1)$,
set $\theta=(\xi\om)^{2-p}$.
There exist 
$c_1,\,\dl\in(0,1)$, depending only on the data $\{\nu, N, p, C_o, C_1\}$, such that if $\xi\om\le 1$ and if
\begin{equation} \label{omi}
| [\pm\left(\mu^\pm-u(\cdot, t_1)\right) \le \xi\om] \cap K_{2\rho}(x_o) | \le c_1(\xi\om)^{\frac{N+2p}{p}} | K_{2\rho}| ,
\end{equation}
then
\begin{equation} \label{akinao}
\pm(\mu^\pm-u) \ge \tfrac14 \xi\om, \quad \text{a.e. in } K_{\frac12\rho}(x_o)\times(t_1+(1-2^{-p})\delta \theta \rho^p, t_1+\dl\theta\rho^p],
\end{equation}
provided $K_{2\rho}(x_o)\times(t_1,t_1+\dl\theta\rho^p]$ is included in $ Q_{R,S}$.
\end{lemma}

\begin{proof}
Let us deal with the case of super-solutions only since the other case is similar. Without loss of generality, we assume $(x_o,t_1)=(0,0)$, as we can always resort to this case by translation. 

We start by defining the sequences, indexed over the set of non-negative integers $\N_0$, 
$$
k_n = \mu^- + \frac{\xi\om}{2} + \frac{\xi\om}{2^{n+1}},
$$
and 
$$
\rho_n = \rho + \frac{\rho}{2^{n}},
$$
and the corresponding cubes and cylinders
$$
K_n=K_{\rho_n} \qquad \mathrm{and} \qquad Q_n=K_n \times ( 0,\delta \theta \rho^p],
$$
where $\delta >0$ will be determined in due course. We also consider the sequence of middle points 
$$
\widetilde{\rho}_n = \frac{\rho_n + \rho_{n+1}}{2}, 
$$
with $\widetilde{K}_n=K_{\widetilde{\rho}_n}$ and $\widetilde{Q}_n=\widetilde{K}_n \times ( 0,\delta \theta \rho^p ]$.

We now write the energy estimate \eqref{Eq:en-est-gen} over $Q_n$, for a cutoff function $\z \in C_o^\infty ( K_n )$ such that $\z(x) \in [0,1]$, $\z \equiv 1$ in $\widetilde{K}_n$ and $| D\z | \leq 2^{n+4}/\rho$, obtaining
\begin{align}
	\essup_{0<t<\delta \theta \rho^p}& \ \int_{\widetilde{K}_n \times\{t\}} (u-k_n)_-^2\,\dx + \iint_{\widetilde{Q}_n}|D(u-k_n)_-|^p\,\dx\dt\label{osaka}\\
	&\le \gm \frac{2^{np}}{\rho^p} \iint_{Q_n} (u-k_n)^{p}_- \,\dx\dt +\int_{K_n\times \{0\}} (u-k_n)_-^2\,\dx\nonumber\\
	&\quad+2\nu\int_{K_n\times\{0\}} (u-k_n)_-\,\dx\nonumber\\
    &\le \gm \frac{2^{np}}{\rho^p}  (\xi\om)^{p}  | [ u < k_n ] \cap Q_n |+ (\xi\om)^2 | [ u (\cdot, 0) < k_n ] \cap K_n | \nonumber\\
	&\quad+2\nu \xi\om | [ u (\cdot, 0) < k_n ] \cap K_n |\nonumber\\
    &\le \gm \frac{2^{np}}{\rho^p}  (\xi\om)^{p}  |[ u < k_n ] \cap Q_n |+ 2 \nu_\ast \xi\om |[ u (\cdot, 0) < k_n ] \cap K_n | \nonumber\\
    &\le \gm \frac{2^{np}}{\rho^p}  (\xi\om)^{p}  | A_n |+ 2^{N+1}\nu_\ast c_1 (\xi\om)^{\frac{N+3p}{p}}  | K_n | \nonumber,
\end{align}
with $\nu_\ast \df= \max \{ 1, \nu \}$.

We used the fact that $\xi\om < 1$ and thus $(\xi\om)^2 < \xi\om$ and, in the last estimate, the following consequence of \eqref{omi}, due to $\rho < \rho_n \leq 2 \rho$:
\begin{align*}
|[u(\cdot, 0) \le k_n ] \cap K_n | & \leq  | [u(\cdot, 0) \le \mu^- + \xi\om] \cap K_{2\rho} | \\
& \leq  c_1 (\xi\om)^{\frac{N+2p}{p}} | K_{2\rho}| \\
& \leq  c_1 (\xi\om)^{\frac{N+2p}{p}} 2^N | K_n |,
\end{align*}
for any $n \in \N_0$. We also denoted
$$A_n \df= [ u < k_n ] \cap Q_n.$$

We now consider an alternative: either, for some $n \in \N_0$,
\begin{equation} \label{samurai}
| A_n | \leq  \frac{c_1(\xi\om)^{\frac{N+p}{p}}}{\delta} | Q_n |
\end{equation}
or we have 
\begin{equation} \label{shogun}
| A_n | >  \frac{c_1(\xi\om)^{\frac{N+p}{p}}}{\delta} | Q_n |,
\end{equation}
for every $n \in \N_0$. If \eqref{samurai} holds, assume $\delta$ is already chosen in terms of the data $\{\nu, N, p, C_o, C_1\}$, and choose $c_1$ as 
\begin{equation} \label{hiromatsu} 
c_1 = c_o 2^{-N-\frac{N}{p}-1}  \delta^{\frac{N}{p}+1},  
\end{equation}
where $c_o$ is determined in Lemma \ref{Lm:DG:1}. Noticing again that $\rho < \rho_n \leq 2 \rho$ and that $k_n > \mu^- + \frac12 \xi\om$, we readily get
$$
| [ u < \mu^- + \tfrac12 \xi\om]  \cap [(0,\delta \theta \rho^p)+Q_{\rho}(\delta \theta)] | \leq | A_n | \leq c_o \delta^{\frac{N}{p}}( \tfrac12 \xi\om )^{\frac{N+p}p}| Q_{\rho} (\delta \theta)|.
$$
Now, Lemma \ref{Lm:DG:1} gives
\begin{equation} \label{toranaga}
u \geq \mu^- + \frac14 \xi\om \qquad \mathrm{in} \ \ K_{\frac12 \rho} \times \left( (1-2^{-p})\delta \theta \rho^p , \delta \theta \rho^p \right].
\end{equation}

If, alternatively, \eqref{shogun} holds, we conclude from the energy estimate \eqref{osaka} that, for every $n \in \N_0$,
\begin{align}
	\essup_{0<t<\delta \theta \rho^p}& \ \int_{\widetilde{K}_n \times\{t\}} (u-k_n)_-^2\,\dx + \iint_{\widetilde{Q}_n}|D(u-k_n)_-|^p\,\dx\dt\label{nagakado}\\
    &\le ( \gm 2^{np} + 2^{N+1}\nu_\ast ) \frac{(\xi\om)^{p}}{\rho^p} | A_n | \nonumber\\
    &\le  \gm 2^{np} \nu_\ast  \frac{(\xi\om)^{p}}{\rho^p} | A_n |, \nonumber
\end{align}
for some $\gm=\gm (N,C_o,C_1) >1$, since \eqref{shogun} implies
$$c_1(\xi\om)^{\frac{N+3p}{p}}   | K_n | < \frac{(\xi\om)^{p}}{\rho^p} | A_n |.$$
We now run De Giorgi's iteration scheme. Take a cutoff function $\phi \in C_o^\infty ( \widetilde{K}_n )$ such that $\phi (x) \in [0,1]$, $\phi  \equiv 1$ in $K_{n+1}$ and $| D\phi  | \leq 2^{n+4}/\rho$, to get, using H\"older's inequality and the parabolic Sobolev inequality (see~\cite[Chapter~I, Proposition~3.1]{dibe-sv}),
\begin{align*}
\frac{\xi\om}{2^{n+2}} | A_{n+1}|  & =   ( k_n-k_{n+1} ) | A_{n+1}|\\
& \leq   \iint_{Q_{n+1}} (u-k_n)_-\,\dx\dt\\
& \leq   \iint_{\widetilde{Q}_n} \phi (u-k_n)_-\,\dx\dt\\
& \leq   \left[ \iint_{\widetilde{Q}_n} | \phi (u-k_n)_- |^{\frac{p(N+2)}{N}}\,\dx\dt \right]^{\frac{N}{p(N+2)}} | A_n |^{1-\frac{N}{p(N+2)}}\\
& \leq   \gm  \left( \iint_{\widetilde{Q}_n} | D [ \phi (u-k_n)_- ] |^{p}\,\dx\dt \right)^{\frac{N}{p(N+2)}} \\
& \quad \cdot \left( \essup_{0<t<\delta \theta \rho^p} \int_{\widetilde{K}_n \times\{ t\}} (u-k_n)_-^{2}\,\dx \right)^{\frac{1}{N+2}} | A_n |^{1-\frac{N}{p(N+2)}}\\
& \leq  \gm  \left( 2^{np} \nu_\ast  \frac{(\xi\om)^{p}}{\rho^p} | A_n | \right)^{\frac{N+p}{p(N+2)}} | A_n |^{1-\frac{N}{p(N+2)}}\\
& =  \gm  \left( 2^{np} \nu_\ast  \frac{(\xi\om)^{p}}{\rho^p} \right)^{\frac{N+p}{p(N+2)}} | A_n |^{1+\frac{1}{N+2}},\\
\end{align*}
for some $\gm=\gm (N,p,C_o,C_1) >1$, where we also used estimate \eqref{nagakado}. We rewrite this inequality for 
$$
Y_n \df= \frac{| A_n |}{| Q_n |},
$$
obtaining the recursive relation
$$
Y_{n+1} \leq  \gm 2^{Npn} \nu_\ast ^{\frac{N+p}{p(N+2)}} \delta^{\frac{1}{N+2}} Y_n^{1+\frac{1}{N+2}}. 
$$
Using the fast geometric convergence of sequences (see~\cite[Chapter~I, Lemma~4.1]{dibe-sv}), we conclude that $Y_n \to 0$ as $n \to \infty$ provided
$$
Y_o \leq \frac{1}{ 2^{Np(N+2)^2}\gm^{N+2} \nu_\ast^{\frac{N+p}{p}} \delta}.   
$$
This is granted by taking 
\begin{equation} \label{mariko}
\delta = \frac{1}{ 2^{Np(N+2)^2}\gm^{N+2} \nu_\ast^{\frac{N+p}{p}}}
\end{equation}
and we obtain 
\begin{equation} \label{fuji}
u \geq \mu^- + \frac12 \xi\om \qquad \mathrm{in} \ \ K_{\rho} \times ( 0 , \delta \theta \rho^p ].
\end{equation}
The result now follows, since no matter which of the alternatives \eqref{samurai} or \eqref{shogun} holds, \eqref{akinao} can be obtained from either \eqref{toranaga} or \eqref{fuji}.

We conclude by tracing the constant dependence: with the choice of $\delta$ in \eqref{mariko}, we obtain from \eqref{hiromatsu}, also taking into account the value of $c_o$ in Remark~\ref{Rmk:c_o},
\[
c_1 =  \frac{\widetilde c_o(N,p,C_o,C_1)}{2^{N+1+N(N+2)^2(N+p)+\frac{N}{p}}\gm^{(N+2)(1+\frac{N}{p})} \nu_\ast ^{(1+\frac{N}{p})^2+(1+\frac{N}{p})}}.
\]
This finishes the proof.
\end{proof}

\section{Sub/super-solutions to the parabolic $p-$Laplacian} \label{S:3}

In this section, we shall consider quasi-linear, parabolic partial differential equations of $p-$Laplacian-type
\begin{equation}\label{Eq:3:1}
\pl_t u-\dvg\bl{\widetilde{A}}(x,t,u, Du) = 0\quad \text{ weakly in }\> E_T,
\end{equation}
where the function $\bl{\widetilde{A}}:E_T\times \rr \times \rr^{N}\to\rn$ satisfies \eqref{Eq:1:2} for $p>2$. The corresponding notion of sub/super-solution is standard, and we refer to \cite[Chapter~II]{dibe-sv}.

Before introducing the main result of this section, we state a helpful lemma that connects the Stefan problem with the parabolic $p-$Laplace equation; its proof can be retrieved from \cite[Lemma~2.3]{BKU-14} or \cite[Lemma~4.3]{Liao-Stefan}. 

\begin{lemma}\label{Lm:sub-solution}
Let $u$ be a local weak solution to the Stefan problem~\eqref{Eq:reg} with structure conditions~\eqref{Eq:1:2}  in $E_T$. 
Then, for any $k>0$, the truncation $(u-k)_+$ is a local weak sub-solution to a parabolic equation of type \eqref{Eq:3:1} in $E_T$. Likewise, for any $k<-\eps$, the truncation $(u-k)_-$ is a local weak sub-solution to a parabolic equation of type \eqref{Eq:3:1} in $E_T$.
\end{lemma}

The main result of this section is the following expansion of positivity. It translates measure-theoretical estimates to pointwise $p-$capacity estimates. The notion of $p-$capacity is discussed in Appendix~\ref{A:1}.

\begin{proposition}\label{Prop:3:1}
Let $k>0$ and $\al\in(0,1)$. Consider a non-negative, local weak super-solution $v$ to the parabolic equation~\eqref{Eq:3:1}, with structure conditions~\eqref{Eq:1:2} in $E_T$. There exist constants $\bar{c},\,\bar{\dl},\,\xi\in(0,1)$ and $b>1$, depending only on the data $\{N,p,C_o,C_1\}$, such that whenever
\[
|[v(\cdot,t)\ge k]\cap K_\rho(x_o)|\ge\al|K_\rho|
\] 
holds true for any $t\in(t_o,t_o+\bar\theta\rho^p]$, where 
$$
\dsty\bar\theta\df=[k \bar\dl \operatorname{cap}_p(K_{\bar{c}\al^2},K_{3})]^{2-p},
$$ 
then, we have
\[
v\ge k \xi\operatorname{cap}_p(K_{\bar{c}\al^2},K_{3}),
\quad\text{a.e. in }\> K_\rho(x_o)\times(t_o+b\xi^{2-p}\bar\theta\rho^p, t_o+2b\xi^{2-p}\bar\theta\rho^p],
\]
provided that  $K_{4\rho}(x_o)\times (t_o,t_o+b\xi^{2-p}\bar\theta(4\rho)^p]\subset E_T$. 
\end{proposition}

\begin{remark}\label{Rmk:xi}\upshape
Once we have Proposition~\ref{Prop:3:1} at our disposal, the constant $\xi$ can be chosen smaller if necessary. As such, a pointwise estimate for $v$ can be actually claimed over longer cylinders as long as they are inside $E_T$.
\end{remark}

\subsection{Preliminary estimates}
Before coming to the actual proof of Proposition~\ref{Prop:3:1}, we present some structural estimates for non-negative super-solutions to \eqref{Eq:3:1} under the structure assumptions \eqref{Eq:1:2}. The following result has appeared in \cite{GSV}, to which we refer for the proof.
\begin{lemma}\label{LBL1}
Consider a non-negative, local weak super-solution $v$ to the parabolic equation~\eqref{Eq:3:1}, with structure conditions~\eqref{Eq:1:2} 
in $E_T$. Assume that $K_{2\rho}(y)\times [\bar t,\bar t+S]\subset E_T$ and that 
$$
\inf_{K_{2\rho}(y)} v(\cdot,\bar t)\geq k,\quad\text{ for some }k>0.
$$ 
Then, there exists $\kappa\in(0,1)$, depending only on the data $\datap$, such that for all $t\in (\bar t,\bar t+S]$ we have
\begin{equation*}
\inf_{K_\rho(y)} v(\cdot,t)\geq \frac{k}{2}\left( 1+\frac{t-\bar t}{\kappa k^{2-p} (2\rho)^p}\right)^\frac{1}{2-p}.
\end{equation*}

\end{lemma}

We also need the following weak Harnack inequality, proved in \cite{Kuusi2008}; see also \cite[Chapter~5, Section~7]{DBGV-mono}.
\begin{theorem}\label{Thm:3:1}
Consider a non-negative, local weak super-solution $v$
to  the parabolic equation~\eqref{Eq:3:1}, with structure conditions~\eqref{Eq:1:2} in $E_T$. Assume that $K_{16\rho}(x_o)\subset E$.  There exist 
constants $c,\,\gm>1$, depending only on the data 
$\datap$, such that, for a.e. $s\in(0,T)$,
\begin{equation}\label{Eq:3:7:1}
\bint_{K_\rho(x_o)}v(x,s)\,\dx\le c
\left(\frac{\rho^p}{T-s}\right)^{\frac1{p-2}}
+\gm\inf_{K_{4\rho}(x_o)}v(\cdot,t),
\end{equation}
for all times 
\begin{equation*}
s+{\txty\frac12}\theta\rho^p\le t\le s+2\theta\rho^p,
\end{equation*}
where 
\begin{equation}\label{Eq:3:7:2}
\theta=\min\left\{c^{2-p}\frac{T-s}{\rho^p}\,,\,\left[
\bint_{K_\rho(x_o)}v(x,s)\,\dx\right]^{2-p}\right\}.
\end{equation}
\end{theorem}
\begin{remark}\label{Rmk:3:2}\upshape
In practice, it will be useful to choose $s$ and $\rho$ to satisfy
\begin{equation}\label{Eq:Harnack:T}
s+ 2c^{p-2}\left[\dsty\bint_{K_\rho(x_o)}v(x,s)\,\dx\right]^{2-p}\rho^p<T,
\end{equation}
that is,
\[
c\left(\frac{\rho^p}{T-s}\right)^{\frac1{p-2}}<\frac1{2^{\frac1{p-2}}}\bint_{K_\rho(x_o)}v(x,s)\,\dx;
\]
then, according to \eqref{Eq:3:7:2}, we have
\begin{align*}
&\theta=\left[\bint_{K_\rho(x_o)}v(x,s)\,\dx\right]^{2-p},
\end{align*}
and also by \eqref{Eq:3:7:1}, 
\begin{equation}\label{WHI}
\bint_{K_\rho(x_o)}v(x,s)\,\dx\le\gm_{\rm H}\inf_{K_{4\rho}(x_o)}v(\cdot,t),
\end{equation}
for all times 
\begin{equation*}
s+{\txty\frac12}\theta\rho^p\le t\le s+2\theta\rho^p.
\end{equation*}
Moreover, $\gm_{\rm H}=\gm/(1-2^{\frac1{2-p}})$, and therefore, the constant is stable as $p\downarrow2$.
\end{remark}

We can show the following estimate by relying on the previous results and working as in the proof of Lemma~3.2 of \cite{GL-deg}.

\begin{lemma}\label{Lm:3:3}
Consider a non-negative, local weak super-solution $v$ 
to the parabolic equation~\eqref{Eq:3:1}, with structure conditions~\eqref{Eq:1:2} in $E_T$. For $k>0$, define
\[
v_k\df=k-(k-v)_+, 
\]
introduce three nested cylinders
\begin{equation*}
\left\{
\begin{array}{ll}
{\bf Q}_2=K_{\rho}(x_o)\times(t_o+\tfrac14\bar\theta\rho^p,t_o+\tfrac34\bar\theta\rho^p],\\[5pt]
{\bf Q}_1=K_{\frac32\rho}(x_o)\times(t_o+\tfrac18\bar\theta\rho^p,t_o+\tfrac78\bar\theta\rho^p],\\[5pt]
{\bf Q}_o=K_{2\rho}(x_o)\times(t_o,t_o+\bar\theta\rho^p],
\end{array}\right.
\end{equation*}
with a parameter $\bar\theta>0$, and assume that ${\bf Q}_o\subset E_T$. For every $\eta\in(0,1)$, there exists a constant $\widetilde{\gm}>1$, depending only on the data $\datap$ and $\eta$, such that
\begin{equation}\label{Eq:3:13}
\begin{aligned}
\bint_{t_o+\frac14\bar\theta\rho^p}^{t_o+\frac34\bar\theta\rho^p}\int_{K_{\frac32\rho}(x_o)}|D(v_k\z)|^p\,\dx\dt\le k  \rho^{N-p} \,{\mathbb I},
\end{aligned}
\end{equation}
where 
\[
{\mathbb I}\df=\widetilde{\gm}\left[\sup_{t_o<t<t_o+\bar\theta\rho^p}\bint_{K_{2\rho}(x_o)} v_k(x,t)\,\dx\right]^{p-1} + \frac{\eta}{\bar\theta^{\frac{p-1}{p-2}}},
\]
and
$\z\in C^1_o({\bf Q}_1;[0,1])$ satisfies $\z=1$ in ${\bf Q}_2$, $|D\z|\le 4/\rho$ and $|\pl_t\z|\le{16}/(\bar\theta\rho^p)$.
\end{lemma}

\begin{proof}
According to \cite[Chapter~3, Lemma~1.1]{DBGV-mono}, $v_k$ is a local weak super-solution to \eqref{Eq:3:1}, with \eqref{Eq:1:2} in $E_T$.
We use the testing function $(k-v_k)\z^p$ in the weak formulation of $v_k$,
modulo a standard Steklov average.
After standard calculations, we get the following energy estimate:
\begin{equation*}
\begin{aligned}
\iint_{{\bf Q}_1}|D(v_k\z)^p|\,\dx\dt&\le \gm k\iint_{{\bf Q}_1}|Dv_k|^{p-1}|D\z|\,\dx\dt+ \gm\iint_{{\bf Q}_1}v_k^p|D\z|^p\,\dx\dt\\
&\quad+ \gm k\iint_{{\bf Q}_1}v_k|\pl_t\z|\,\dx\dt,
\end{aligned}
\end{equation*}
for some $\gm=\gm(p,C_o,C_1)$.
This is precisely what one obtains in the proof of \cite[Lemma 3.2]{GL-deg}. Starting from this, one can reproduce the arguments and arrive at
\[
\iint_{{\bf Q}_1}|D(v_k\z)|^p\,\dx\dt\le C_\eta k\bar{\theta}\rho^N\left[\sup_{t_o<t<t_o+\bar\theta\rho^p}\bint_{K_{2\rho}(x_o)} v_k(x,t)\,\dx\right]^{p-1} + k\bar{\theta}\rho^N\frac{\eta}{\bar\theta^{\frac{p-1}{p-2}}},
\]
for any $\eta\in(0,\frac12)$ and some $C_{\eta}>1$.
A consequence of the previous estimate is
\begin{align*}
\bint_{t_o+\frac14\bar\theta\rho^p}^{t_o+\frac34\bar\theta\rho^p}&\int_{K_{\frac32\rho}(x_o)}|D(v_k\z)|^p\,\dx\dt\\
&\le\frac1{\frac12\bar\theta\rho^p}\int_{t_o+\frac18\bar\theta\rho^p}^{t_o+\frac78\bar\theta\rho^p}\int_{K_{\frac32\rho}(x_o)}|D(v_k\z)|^p\,\dx\dt\\ 
&\le 2C_\eta k\rho^{N-p}\left[\sup_{t_o<t<t_o+\bar\theta\rho^p}\bint_{K_{2\rho}(x_o)} v_k(x,t)\,\dx\right]^{p-1} + k\rho^{N-p}\frac{2\eta}{\bar\theta^{\frac{p-1}{p-2}}}.
\end{align*}
The proof is concluded by redefining $2\eta$ as $\eta$ and choosing $\widetilde\gm=2C_{\eta}$.
\end{proof}

\begin{remark}\label{Rmk:3:3}\upshape
The role of $\bar\theta$ becomes significant when it restores the homogeneity of the estimate~\eqref{Eq:3:13}.
In fact, if we let
\begin{equation*}
\bar\theta=(\dl k)^{2-p},
\end{equation*}
for some $\dl\in(0,1)$ to be determined later,
 it is straightforward to see that
\begin{equation*}
\frac1{\bar\theta^{\frac{p-1}{p-2}}}=(\dl k)^{p-1}<k^{p-1}.
\end{equation*}
Moreover, by definition of $v_k$, we  have
\[
\sup_{t_o<t<t_o+\bar\theta\rho^p}\bint_{K_{2\rho}(x_o)} v_k(x,t)\,\dx\le k.
\]
Therefore, under the choice of $\bar\theta$, we have ${\mathbb I}\le 2\widetilde{\gm} k^{p-1}$. Keep in mind that $\widetilde{\gm}=\widetilde{\gm}(\eta)$.
\end{remark}

\subsection{Proof of Proposition~\ref{Prop:3:1}}
For some $\varep\in(0,\frac12)$ to be determined later, without loss of generality, we can assume that $\frac12\varep^{-p}$ is an integer; hence, we divide the interval $(t_o+\frac14\bar\theta\rho^p,t_o+\frac34\bar\theta\rho^p]$ into $\frac12 \varep^{-p}$ sub-intervals, each of length $\bar\theta(\varep\rho)^p$. In the present subsection, we let $\bar\theta=(\dl k)^{2-p}$, for some $\dl\in(0,1)$, so that Remark~\ref{Rmk:3:3} applies. The actual form of $\dl$ will become clear in the course of the proof. 

By \eqref{Eq:3:13}, there must exist one of these sub-intervals, say $(t_*,t_*+\bar\theta(\varep\rho)^p]$, such that
\begin{equation}\label{Eq:3:14+}
\bint_{t_*}^{t_*+\bar\theta(\varep\rho)^p}\int_{K_{\frac32\rho}(x_o)}|D(v_k\z)|^p\,\dx\dt\le  k \rho^{N-p}\,{\mathbb I},
\end{equation}
which immediately implies the following estimate:
\begin{equation}\label{Eq:3:14}
\begin{aligned}
&\bint_{t_*}^{t_*+\frac12\bar\theta(\varep\rho)^p}\int_{K_{\frac32\rho}(x_o)}|D(v_k\z)|^p\,\dx\dt\le2  k \rho^{N-p}\,{\mathbb I}.
\end{aligned}
\end{equation}
Note that these estimates do not require any quantitative knowledge of $\varep$, which is still to be selected. From \eqref{Eq:3:14}, recalling that $\z\equiv1$ in 
${\bf Q}_2=K_{\rho}(x_o)\times(t_o+\frac14\bar\theta\rho^p,t_o+\frac34\bar\theta\rho^p]$, and also that ${\mathbb I}\le 2\widetilde\gm k^{p-1}$, $\widetilde{\gm}=\widetilde{\gm}(\eta)$ from Remark~\ref{Rmk:3:3}, we infer that there exists
some 
$$\bar t\in[t_*,t_*+\tfrac12\bar\theta(\varep\rho)^p]$$ 
satisfying
\[
\int_{K_{\rho}(x_o)\times\{\bar t\}}|Dv_k|^p\,\dx\dt\le4\widetilde\gm k^p \rho^{N-p},
\]
while the measure-theoretical information yields
\begin{equation*}
|[v_k(\cdot,\bar t)=k]\cap K_\rho(x_o)|\ge\al|K_\rho|.
\end{equation*}
The last two estimates permit us to apply Lemma~\ref{Lm:0:3:1}.
Indeed, using the lemma for $v_k(\cdot, \bar{t})$, with $\dl=\lm=\frac12$, we find
\begin{equation}\label{Eq:eps:0}
\varep=\bar{c}\al^2
\end{equation} as in \eqref{Eq:0:3:7}, with $\bar{c}\in(0,\frac12)$ depending only on the data $\datap$ and $\eta$,
such that
\begin{equation}\label{Eq:K-eps-meas}
|[v_k(\cdot,\bar t)>\tfrac12 k]\cap K_{\varep\rho}(y)|>\tfrac12|K_{\varep\rho}|,
\end{equation}
for some $K_{\varep\rho}(y)\subset K_\rho(x_o)$. 

Since $v_k$ is a weak super-solution to \eqref{Eq:3:1}, with structure conditions \eqref{Eq:1:2} in $E_T$, we apply the weak Harnack inequality \eqref{WHI} to $v_k$ on $K_{\varep\rho}(y)$ at time $\bar t$. As a result, we obtain
\begin{equation}\label{Eq:3:15}
\bint_{K_{\varep\rho}(y)} v_k(x,\bar t)\,\dx\le \gm_{{\rm H}}\inf_{K_{4\varep\rho}(y)}v_k(\cdot,t),
\end{equation}
for all times 
\[
\bar t+\tfrac12 \theta(\varep\rho)^p\le t \le\bar t+2 \theta(\varep\rho)^p,\quad\text{ where }\quad \theta=\left[\bint_{K_{\varep\rho}(y)}v_k(x,\bar t)\,\dx\right]^{2-p},
\]
provided that \eqref{Eq:Harnack:T} holds, that is, $\bar{t}+2c^{p-2}\theta(\varep\rho)^{p}<T$ in the current set-up. 

Let us take a closer look at the above conclusion using  \eqref{Eq:K-eps-meas}. 
Notice first that
\begin{align*}
\bint_{K_{\varep\rho}(y)} v_k(x,\bar t)\,\dx
\ge\frac1{|K_{\varep\rho}|}\int_{K_{\varep\rho}(y)\cap[v_k>\frac12 k]} v_k(x,\bar t)\,\dx>\tfrac14 k,
\end{align*}
and, on the other hand, by definition of $v_k$, we have
\[
\bint_{K_{\varep\rho}(y)}v_k(x,\bar t)\,\dx\le k.
\]
Hence, we can estimate $\theta$ as
\begin{equation}\label{Eq:Est:und-thet}
k^{2-p}\le \theta=\left[\bint_{K_{\varep\rho}(y)} v_k(x,\bar t)\,\dx\right]^{2-p}<(\tfrac14 k)^{2-p},
\end{equation}
and, from \eqref{Eq:3:15}, we infer that
\begin{equation}\label{Eq:Harnack:0}
v_k\ge\frac{k}{4\gm_{\rm H}}\df{=}\bar\eta k, \quad\text{ a.e. in } K_{\varep\rho}(y)\times[\bar t+\tfrac12 \theta(\varep\rho)^p,\bar t+2 \theta(\varep\rho)^p].
\end{equation}
By \eqref{Eq:Est:und-thet}, $\bar{t}\le t_o+\bar\theta\rho^p$, and $\bar\theta=(\dl k)^{2-p}$, the requirement $\bar{t}+2c^{p-2}\theta(\varep\rho)^{p}<T$ is fulfilled if 
\begin{equation}\label{Eq:Harnack:T:0}
t_o+2^p c^{p-2} \bar{\theta}\rho^p<T.
\end{equation}

No matter the value of $\bar t\in[t_*,t_*+\frac12\bar\theta(\varep\rho)^p]$, provided $\dl$ in $\theta=(\dl k)^{2-p}$ is taken sufficiently small, we can surely conclude that 
\begin{equation}\label{Eq:time-incl}
K_{\varep\rho}(y)\times[\bar t+\tfrac12 \theta(\varep\rho)^p,\bar t+2 \theta(\varep\rho)^p]\subset K_\rho(x_o)\times(t_*,t_*+\bar\theta(\varep\rho)^p].
\end{equation}
Indeed, since
\[
\bar t+2 \theta(\varep\rho)^p\le t_*+\tfrac12\bar\theta(\varep\rho)^p+2 \theta(\varep\rho)^p,
\]
in order to have the wanted inclusion, it suffices to require
\begin{align*}
t_*+\tfrac12\bar\theta(\varep\rho)^p+2 \theta(\varep\rho)^p\le t_*+\bar\theta(\varep\rho)^p,\quad\textit{i.e.},\quad
 \theta\le\tfrac14\bar\theta,
\end{align*}
which is implied if
\begin{equation}\label{Eq:up-bd-dl}
\dl\le 4^{\frac{p-1}{2-p}},
\end{equation}
using $\bar\theta=(\dl k)^{2-p}$ and the right-hand side of \eqref{Eq:Est:und-thet}.
This sets a smallness requirement on $\dl$, while a more precise choice of $\dl$ is yet to come.

The inclusion \eqref{Eq:time-incl} allows us to estimate
\[
\bint_{t_*}^{t_*+\bar\theta(\varep\rho)^p}\int_{K_{\frac32\rho}(x_o)}|D(v_k\z)|^p\,\dx\dt\ge\frac{\frac32 \theta(\varep\rho)^p}{\bar\theta(\varep\rho)^p}\,\,\bint_{\bar t+\frac12 \theta(\varep\rho)^p}^{\bar t+2 \theta(\varep\rho)^p}\int_{K_{\frac32\rho}(x_o)}|D(v_k\z)|^p\,\dx\dt,
\]
and hence, by \eqref{Eq:3:14+}, we have
\[
\frac{3 \theta}{2\bar\theta}\,\,\bint_{\bar t+\frac12 \theta(\varep\rho)^p}^{\bar t+2 \theta(\varep\rho)^p}\int_{K_{\frac32\rho}(x_o)}|D(v_k\z)|^p\,\dx\dt\le  k \rho^{N-p}\,{\mathbb I}.
\]
The last display, together with \eqref{Eq:Harnack:0}, gives some $\widetilde t\in(\bar t+\tfrac12 \theta(\varep\rho)^p,\bar t+2 \theta(\varep\rho)^p]$, such that
\[
\frac{3 \theta}{2\bar\theta}\,\int_{K_{\frac32\rho}(x_o)\times\{\widetilde t\}}|D(v_k\z)|^p\,\dx\le  k \rho^{N-p}\,{\mathbb I}\quad\text{ and }\quad \inf_{K_{\varep\rho}(y)}v_k(\cdot,\widetilde{t})\ge\bar\eta k.
\]
By properties of the $p-$capacity (see~Appendix~\ref{A:1}), $K_{\varep\rho}(y)\subset K_{\rho}(x_o)$, and $K_{\frac32\rho}(x_o)\subset K_{3\rho}(y)$, we obtain
\[
  k \rho^{N-p}\,{\mathbb I}\ge \frac{3 \theta}{2\bar\theta}(\bar\eta k)^p\operatorname{cap}_p(K_{\varep\rho}(y),K_{\frac32\rho}(x_o))\ge \frac{3 \theta}{2\bar\theta}(\bar\eta k)^p\operatorname{cap}_p(K_{\varep\rho}(y),K_{3\rho}(y)).
\]
By $\bar\theta=(\dl k)^{2-p}$ and the estimate of $\theta$ in \eqref{Eq:Est:und-thet}, we further compute
\begin{align*}
  k \rho^{N-p}\,{\mathbb I}&\ge \frac{3 \dl^{p-2}}{2}(\bar\eta k)^p\operatorname{cap}_p(K_{\varep\rho}(y),K_{3\rho}(y)),\quad \textit{i.e.},\\
{\mathbb I}&\ge\frac{3 \bar\eta^p}{2 } \,\dl^{p-2}\frac{\operatorname{cap}_p(K_{\varep\rho}(y),K_{3\rho}(y))}{\rho^{N-p}} k^{p-1}.
\end{align*}
Recalling the definition of ${\mathbb I}$ and also the scaling property \eqref{Eq:cap-scale} of $p-$capacity, the last inequality yields 
\[
\frac{3 \bar\eta^p}{2 } \,\dl^{p-2} \operatorname{cap}_p(K_{\varep},K_{3}) \, k^{p-1}\le\widetilde\gm \left[\sup_{t_o<t<t_o+\bar\theta\rho^p}\bint_{K_{2\rho}(x_o)} v_k(x,t)\,\dx\right]^{p-1} +\eta \dl^{p-1} k^{p-1}.
\]
At this point, we select the parameters $\dl$ and $\eta$ so that the second term on the right-hand side will be absorbed into the left. This leads to the choices
\begin{equation}\label{Eq:dl-eta}
\dl(\varep)=\bar\dl\operatorname{cap}_p(K_{\varep},K_{3}) ,\qquad \eta=\frac{3 \bar\eta^p}{4\bar{\dl}}.
\end{equation}
Here, $\bar\dl\in(0,1)$ can be selected in terms of $\{N, p\}$ only to guarantee the previous smallness requirement \eqref{Eq:up-bd-dl} on $\dl$. In fact, by \eqref{Eq:cap-bd}, we have $\operatorname{cap}_p(K_{\varep},K_{3})\le \gm(N,p)$, and hence it suffices to take
\[ 
\bar\dl=\frac{4^{\frac{p-1}{2-p}}}{\gm},
\] 
whereas $\bar\eta=1/(4\gm_{\rm H})$ is from \eqref{Eq:Harnack:0}. Thus, at this point, $\eta$ has been determined in \eqref{Eq:dl-eta} by the data $\datap$, and also $\varep$ in \eqref{Eq:eps:0}.
Consequently, these choices lead to
\begin{align*}
\frac{3 \bar\eta^p \bar{\dl}^{p-2}}{4\widetilde\gm} \,[\operatorname{cap}_p(K_{\varep},K_{3})]^{p-1} k^{p-1}&\le \left[\sup_{t_o<t<t_o+\bar\theta\rho^p}\bint_{K_{2\rho}(x_o)} v_k(x,t)\,\dx\right]^{p-1},
\end{align*}
\textit{i.e.},
\begin{align*}
\frac{\eta \bar{\dl}^{p-1}}{\widetilde\gm} \,[\operatorname{cap}_p(K_{\varep},K_{3})]^{p-1} k^{p-1}&\le \left[\sup_{t_o<t<t_o+\bar\theta\rho^p}\bint_{K_{2\rho}(x_o)} v_k(x,t)\,\dx\right]^{p-1};
\end{align*}
in other words, we have
\begin{equation}\label{Eq:low-bd-1}
\eta_* \bar{\dl} \operatorname{cap}_p(K_{\varep},K_{3})\, k \le \sup_{t_o<t<t_o+\bar\theta\rho^p}\bint_{K_{2\rho}(x_o)} v_k(x,t)\,\dx,
\end{equation}
where 
\begin{equation}\label{Eq:low-bd-2}
\bar\theta=(\dl k)^{2-p}=[\operatorname{cap}_p(K_{\varep},K_{3})]^{2-p} (\bar\dl k)^{2-p}\qquad\text{ and }\qquad \eta_*= \left(\frac{\eta}{\widetilde\gm}\right)^{\frac1{p-1}}.
\end{equation}
Note that $\bar{\eta}=1/(4\gm_{\rm H})$ from \eqref{Eq:Harnack:0}, $\eta$ is given by \eqref{Eq:dl-eta}, and $\widetilde\gm=\widetilde\gm(\eta)$ is fixed once $\eta$ is chosen in \eqref{Eq:dl-eta}. Hence,  $\eta_*\in(0,1)$ depends only on the data  $\datap$.

Now, we analyze the consequences of   \eqref{Eq:low-bd-1}, together with the weak Harnack inequality. Recall the choice of $\bar\theta$ in \eqref{Eq:low-bd-2}, and for simplicity, denote
\[
\widetilde{\dl}(\varep)\df{=}\operatorname{cap}_p(K_{\varep},K_{3})\quad\Rightarrow\quad\dl(\varep)=\bar{\dl}\widetilde{\dl}(\varep),\quad\bar\theta=[\bar{\dl}\widetilde{\dl}(\varep) k]^{2-p}.
\]
Let $t_1\in[t_o,t_o+\bar\theta\rho^p]$ be the level where the supremum in \eqref{Eq:low-bd-1} is attained, so that 
\begin{equation}\label{Eq:Harnack:1}
\eta_*\dl(\varep) k\le\bint_{K_{2\rho}(x_o)}v_k(x,t_1)\,\dx.
\end{equation}
The weak Harnack inequality \eqref{WHI} applied to $v_k$ on $K_{2\rho}(x_o)$ at time $t_1$ yields
\begin{equation}\label{Eq:3:15bis}
\bint_{K_{2\rho}(x_o)} v_k(x,t_1)\,\dx\le \gm_{\rm H}\inf_{K_{4\rho}(x_o)}v_k(\cdot,t),
\end{equation}
for all times 
\[
t_1+\tfrac12\theta_1(2\rho)^p \leq t\le t_1+2\theta_1(2\rho)^p,\quad\text{ where }\quad\theta_1=\left[\bint_{K_{2\rho}(x_o)}v_k(x,t_1)\,\dx\right]^{2-p},
\]
provided that \eqref{Eq:Harnack:T} holds, \textit{i.e.}, in the current set-up,
\begin{equation}\label{Eq:Harnack:T:1}
t_1 + 2 c^{p-2}\theta_1(2\rho)^p<T.
\end{equation}
Assume that \eqref{Eq:Harnack:T:1} holds for the moment; then,
combining \eqref{Eq:Harnack:1} and \eqref{Eq:3:15bis}, we obtain
\begin{equation}\label{Eq:Harnack:2}
\inf_{K_{4\rho}(x_o)} v_k(\cdot,\widehat t)\ge\frac{\eta_*}{\gm_{\rm H}}\dl(\varep)k,\quad\text{ where }\ \widehat t\df=t_1+2\theta_1(2\rho)^p.
\end{equation}

Now, let us examine what is required for \eqref{Eq:Harnack:T:1} to hold true.
From \eqref{Eq:Harnack:1}, it is apparent that 
\[
\theta_1=\left[\bint_{K_{2\rho}(x_o)}v_k(x,t_1)\,\dx\right]^{2-p}\le(\eta_*\dl(\varep) k)^{2-p}=\eta_*^{2-p}\bar\theta.
\]
Recalling that $t_1\in [t_o,t_o+\bar\theta\rho^p]$, we estimate
\[
t_1 + 2 c^{p-2}\theta_1(2\rho)^p\le t_o + \bar\theta\rho^p +2 c^{p-2}\eta_*^{2-p}\bar\theta(2\rho)^p=t_o+(1+2^{p+1} c^{p-2}\eta_*^{2-p})\bar\theta\rho^p.
\]
Hence, condition \eqref{Eq:Harnack:T:1} amounts to requiring
\begin{equation}\label{Eq:Harnack:T:2} 
t_o+(1+2^{p+1} c^{p-2}\eta_*^{2-p})\bar\theta\rho^p<T,
\end{equation}
which we may assume.
Likewise, we also have
\begin{equation}\label{Eq:t-hat}
\widehat{t}=t_1+2\theta_1(2\rho)^p\le t_o+\bar\theta\rho^p+2^{p+1}\eta_*^{2-p}\bar\theta\rho^p=t_o+(1+2^{p+1}\eta_*^{2-p})\bar\theta\rho^p.
\end{equation}

Finally, we aim to use \eqref{Eq:Harnack:2} and apply Lemma~\ref{LBL1}. Indeed, consider $t\in(t_o+4^p\eta_*^{2-p}\bar\theta\rho^p,t_o+8^p\eta_*^{2-p}\bar\theta\rho^p]$. Then, estimate \eqref{Eq:t-hat} gives that $0<t-\widehat t\le8^p\eta_*^{2-p}\bar\theta\rho^p$, and also that
\[
1+\frac{t-\widehat t}{\kappa[\frac{\eta_*}{\gm_{\rm H}}\dl(\varep)k]^{2-p}(4\rho)^p}\le 1+ \frac{2^p}{\kappa \gm_{\rm H}^{p-2}}.
\]
Hence, a straightforward application of Lemma~\ref{LBL1} yields 
\begin{equation}\label{Eq:Harnack:3}
\inf_{K_\rho(x_o)} v_k(\cdot,t)\ge\frac{\eta_*\dl(\varep) k}{2\gm_{\rm H}} \left(1+ \frac{2^p}{\kappa \gm_{\rm H}^{p-2}}\right)^{\frac1{2-p}}\df{=}\widetilde{\eta} [\eta_*\dl(\varep)k],
\end{equation}
for all $t\in(t_o+4^p[\eta_*\dl(\varep)k]^{2-p}\rho^p,t_o+8^p[\eta_*\dl(\varep)k]^{2-p}\rho^p]$,
provided that
\begin{equation}\label{Eq:Harnack:T:3}
t_o+8^p\eta_*^{2-p}\bar\theta\rho^p<T.
\end{equation}
Now, it suffices to define $\xi=\widetilde{\eta}\eta_*$ and choose $b=\max\{4^p\widetilde{\eta}^{p-2},\,4^p(\widetilde{\eta} c)^{p-2}\}$. In this way the desired pointwise estimate is given in \eqref{Eq:Harnack:3} if we require $t_o+ b\xi^{2-p}\bar\theta(4\rho)^p<T$, so that the conditions \eqref{Eq:Harnack:T:0}, \eqref{Eq:Harnack:T:1}, \eqref{Eq:Harnack:T:2} and \eqref{Eq:Harnack:T:3} are satisfied. Taking into account the definitions of $\varep$ in \eqref{Eq:eps:0} and $\dl(\varep)$ in \eqref{Eq:dl-eta}, we conclude.

\section{Proof of Theorem~\ref{Thm:1}}\label{S:proof-thm}

Let $u_{\eps}$ be a solution to the regularized problem \eqref{Eq:reg}; we recall that $u_{\eps}$ is locally H\"older continuous, hence, in particular, defined everywhere. From now on, we will omit the $\eps$ for simplicity. Moreover, we will focus primarily on the case \underline{$p=N$}, discussing the changes needed to cover the remaining cases at the end of the section. 

Let $b>1$, $\bar c,\,\bar\delta,\,\xi\in(0,1)$ be the constants stipulated in Proposition~\ref{Prop:3:1}; they depend only on $\{N,C_o,C_1\}$, whereas the constant  $c_1\in(0,1)$ postulated in Lemma~\ref{Lm:DG-new}  depends on $\{\nu,N,C_o,C_1\}$. 
Suppose there are parameters $\rho>0$ and $ \om\in(0,1)$ such that
\begin{equation}\label{osc-control}
Q_{8\rho}(b\xi^{2-N}\bar\theta)\subset E_T,\qquad\osc_{Q_{8\rho}(b\xi^{2-N}\bar\theta)} u\le\om
\end{equation}
where
\begin{equation}\label{al-theta}
\bar\theta\df=\left[\left(\tfrac18\om\right)\bar\delta\,{\rm cap}_N(K_{\bar c\al^2},K_3)\right]^{2-N},\qquad\al\df=c_1\left(\tfrac18\om\right)^3.
\end{equation}
Here and in what follows, we assume up to a translation that $(x_o,t_o)=(0,0)$.
Moreover, we consider the quantities
\[
\mu^+=\sup_{Q_{8\rho}(b\xi^{2-N}\bar\theta)} u,\qquad\mu^-=\inf_{Q_{8\rho}(b\xi^{2-N}\bar\theta)} u,
\]
and assume that
\[
\mu^+-\mu^-\ge\tfrac12\om,
\]
dealing with the opposite situation at the end. Due to this assumption, we have that one of the following \textbf{two cases} must hold: either 
\begin{equation}\label{I-alt}
\mu^+-\tfrac18\om\ge\tfrac18\om
\end{equation}
or
\begin{equation}\label{II-alt}
\mu^-+\tfrac18\om\le-\tfrac18\om.
\end{equation}

Let us first suppose that the \textbf{first case} \eqref{I-alt} holds true. We perform the reduction of oscillation in the following two subsections under the first case. For that, we let
\[
I\df=(-2b\xi^{2-N}\bar\theta\rho^N,-b\xi^{2-N}\bar\theta\rho^N],
\]
and further consider the following \textbf{two alternatives}:
\begin{align}
\exists\,\bar t\in I\quad\text{ such that }\quad& | [ u(\cdot,\bar t)\le \mu^-+\tfrac18\om ] \cap K_\rho | \le c_1(\tfrac18\om)^3|K_\rho|,\label{I-alt-i}\\
&\nonumber\\
\forall\,t\in I\quad\text{ we have }\quad& | [ u(\cdot,t)\le\mu^-+\tfrac18\om ] \cap K_\rho | > \underbrace{c_1(\tfrac18\om)^3}_{=\al}|K_\rho|.\label{I-alt-ii}
\end{align}

\subsection{First case -- first alternative}\label{SS:4:1}
In this section, we assume that the \textbf{first case} \eqref{I-alt} holds true and consider the \textbf{first alternative} \eqref{I-alt-i}. Then, Lemma~\ref{Lm:DG-new} for super-solutions yields
\[
u\big(\cdot,\bar t+\tfrac{\delta}{2^N}(\tfrac18\om)^{2-N}\rho^N\big)\ge\mu^-+\tfrac1{32}\om \quad\text{ in }\,K_{\frac14\rho},
\]
where $\delta\in(0,1)$ is the quantity stipulated in Lemma~\ref{Lm:DG-new}, which depends only on the data $\{\nu,N,C_o,C_1\}$. Then, for any $\widetilde\xi\in(0,\tfrac1{32}]$ we obviously have
\[
u\big(\cdot,\bar t+\tfrac{\delta}{2^N}(\tfrac18\om)^{2-N}\rho^N\big)\ge\mu^-+\widetilde\xi\om \quad\text{ in }\,K_{\frac14\rho}.
\]
We want to apply Lemma~\ref{Lm:DG:initial:1} for super-solutions, relying on such pointwise information with time level $t_1=\bar t+\frac{\delta}{2^N}(\frac18\om)^{2-N}\rho^N$ for a proper $\widetilde\xi$ to be chosen; provided we let
$\theta=\gamma_o(\widetilde\xi\om)^{2-N}$ in Lemma~\ref{Lm:DG:initial:1},
the measure-theoretical condition of Lemma~\ref{Lm:DG:initial:1} is automatically satisfied, and we conclude that
\begin{equation}\label{Eq:lower-bd:1}
    u\ge\mu^-+\tfrac12\widetilde\xi\om
\end{equation}
in 
\begin{equation}\label{Eq:cylinder:1}
K_{\frac1{8}\rho}
\times\Big( \bar t+\tfrac\delta{2^N}(\tfrac18\om)^{2-N}\rho^N,\bar t+\tfrac\delta{2^N}(\tfrac18\om)^{2-N}\rho^N+\gamma_o(\widetilde\xi\om)^{2-N}\tfrac1{4^N}\rho^N\Big],
\end{equation}
where $\widetilde\xi$ is still to be determined.

Considering the location of $\bar t$ in $I$, and taking into account the definition of $\bar\theta$ in \eqref{al-theta}, in order for the previous estimate \eqref{Eq:lower-bd:1} to reach $t=0$, it suffices to choose $\widetilde\xi$ such that 
\begin{align}\label{Eq:t-bar:1}
    \tfrac\delta{2^N}(\tfrac18\om)^{2-N}\rho^N+\gamma_o(\widetilde\xi\om)^{2-N}\tfrac1{4^N}\rho^N&\ge 2b\xi^{2-N}\bar\theta \rho^N\\ \nonumber
    &=
2b\xi^{2-N}\big[(\tfrac18\om)\bar\delta\,{\rm cap}_N(K_{\bar c\al^2},K_3)\big]^{2-N}\rho^N,
\end{align}
which is the case if, discarding the term containing $\dl$,
\[
\widetilde\xi\le \Big(\frac{\gamma_o}{b}\Big)^{\frac1{N-2}}2^{\frac{5(1-N)}{N-2}} \xi\, \bar\delta\, [{\rm cap}_N(K_{\bar c\al^2},K_3) ].
\]
By the capacity estimate of \eqref{Eq:cap-dim:1}$_{2}$, and the definition of $\al$ in \eqref{al-theta}, we calculate 
\begin{equation}\label{Eq:cap-ln}
    {\rm cap}_N(K_{\bar c\al^2},K_3)=b_2\big|\ln\big[\bar c c_1^2(\tfrac18\om)^6\big]\big|^{-(N-1)}
=6^{-(N-1)}\, b_2|\ln(\bar\xi\om)|^{-(N-1)},
\end{equation}
for some $b_2$ depending only on $N$, where we denoted 
\[
\bar\xi^6\df=\frac{\bar c c_1^2}{8^6}
\]
for simplicity. Note that $\bar\xi=\gm( N, C_o, C_1)/\nu_*^2$ because of the dependencies of $c_1$.

It is apparent that $\dsty\bar\xi\in(0,1)$. Then, we can choose
\begin{equation}\label{choice-xi}
  \widetilde\xi=\underbrace{\Big(\frac{\gamma_o}{b}\Big)^{\frac1{N-2}}2^{\frac{5(1-N)}{N-2}} \xi\, \bar\delta\, 6^{-(N-1)}\, b_2}_{\df= 2\bar\eta}|\ln(\bar\xi\om)|^{-(N-1)}.  
\end{equation}
Note that $\bar\eta=\bar\eta( N, C_o, C_1)$ is independent of $\nu$. Next, by imposing  $\xi<2^{-18}$ (see~Remark~\ref{Rmk:xi}), we further estimate that
\begin{align*}
  b\xi^{2-N}\bar\theta\rho^N&=b\xi^{2-N}\big[(\tfrac18\om)\bar\delta\,{\rm cap}_N(K_{\bar c\al^2},K_3)\big]^{2-N}\rho^N\\
  &\overset{\eqref{Eq:cap-ln}}{=}b\xi^{2-N}\Big[(\tfrac18\om)\bar\delta\, 6^{-(N-1)}\, b_2|\ln(\bar\xi\om)|^{-(N-1)}\Big]^{2-N}\rho^N\\
  &>(32^N+1)\theta(\tfrac18\rho)^N , 
\end{align*}

\noindent where we defined
\begin{equation}\label{Eq:theta}
\theta\df=\bigg(\frac{\om}{|\ln(\bar\xi\om)|^{N-1}}\bigg)^{2-N}.
\end{equation}
As a result of the last estimate, we have
\[
-b\xi^{2-N}\bar\theta\rho^N+\tfrac\delta{2^N}(\tfrac18\om)^{2-N}\rho^N<-\theta(\tfrac18\rho)^N.
\]
The above inequality, jointly with \eqref{Eq:t-bar:1}, guarantees that the cylinder \eqref{Eq:cylinder:1} includes $Q_{\frac1{8}\rho}(\theta)$, no matter where $\bar t$ is in the interval $I$. Therefore, the pointwise estimate \eqref{Eq:lower-bd:1} can be claimed in $Q_{\frac1{8}\rho}(\theta)$, and this, in turn, yields the oscillation estimate
\begin{equation}\label{osc-1}
\osc_{Q_{\frac1{8}\rho}(\theta)}u\le(1-\eta)\om,
\end{equation}
where $\theta$ is defined in \eqref{Eq:theta} and
\begin{equation}\label{theta-eta}
\begin{aligned}
\eta\df=\frac{\bar\eta}{\left|\ln(\bar\xi\om)\right|^{N-1}}.
\end{aligned}
\end{equation} 
Notice that $\bar\eta$ depends only on the data $\{N,C_o,C_1\}$ while $\bar\xi$ depends additionally on $\nu$.

\subsection{First case -- second alternative}\label{SS:4:2}
In this subsection, we still assume that the \textbf{first case} \eqref{I-alt} holds true. However, we consider the \textbf{second alternative} \eqref{I-alt-ii}. Since $\mu^+-\tfrac18\om\ge\mu^-+\tfrac18\om$, we can rephrase it in the following way:
\[
\forall\,t\in I\quad\text{ we have }\quad|[u(\cdot,t)\le\mu^+-\tfrac18\om]\cap K_\rho|> c_1(\tfrac18\om)^3|K_\rho|.
\]
We introduce the new function
\[
v\df=\tfrac18\om-(u-k)_+ , \quad\text{ with }\,\,k=\mu^+-\tfrac18\om.
\]
Since \eqref{I-alt} ensures that $\mu^+-\frac18\om>0$,
by Lemma~\ref{Lm:sub-solution}, the above-defined $v$ is a non-negative, weak super-solution to the parabolic equation~\eqref{Eq:3:1}, with structure conditions~\eqref{Eq:1:2} in $Q_{8\rho}(b\xi^{2-N}\bar\theta)$. Moreover, in terms of $v$, the previous measure-theoretical information can be rephrased as
\[
\forall\,t\in I,\qquad |[v(\cdot,t)=\tfrac18\om]\cap K_\rho|> \underbrace{c_1(\tfrac18\om)^3}_{=\al}|K_\rho|.
\]
In particular, such information holds for any $t\in J$, where 
\[
J\df=(-2b\xi^{2-N}\bar\theta\rho^N,-(2b\xi^{2-N}-1)\bar\theta\rho^N],
\]
and $\bar\theta$ is defined in \eqref{al-theta}. We can then apply Proposition~\ref{Prop:3:1}, with $t_o=-2b\xi^{2-N}\bar\theta\rho^N$ and $k=\frac18\om$, and conclude that 
\[
v\ge\tfrac18\om\,\xi\,{\rm cap}_N(K_{\bar c\al^2},K_3)\overset{\eqref{Eq:cap-ln}}{=}\tfrac18\om\,\xi\,6^{-(N-1)}\,b_2|\ln(\bar\xi\om)|^{-(N-1)}
\]
in $K_{\frac12\rho}\times(-b\xi^{2-N}\bar\theta\rho^N,0]$. Taking into account the definition of $v$, this yields the following reduction of oscillation:
\begin{equation}\label{osc-2}
    \osc_{Q_{\frac1{2}\rho}(\theta)}u\le(1-\eta)\om,
\end{equation}
where $\theta$ is as in \eqref{Eq:theta} and
\begin{equation}\label{theta-eta-bis}
\begin{aligned}
\eta\df=\underbrace{\tfrac18 \xi\, 6^{-(N-1)}\, b_2}_{\df=\bar\eta}|\ln(\bar\xi\om)|^{-(N-1)} = \frac{\bar\eta}{|\ln(\bar\xi\om)|^{N-1}}.
\end{aligned}
\end{equation}
Notice that the above-defined $\bar\eta$ depends only on the data $\{N,C_o,C_1\}$ and might differ from the one in \eqref{theta-eta}.
However, we tacitly take the smaller of the two.
Therefore, we conclude that when the \textbf{first case} \eqref{I-alt} holds, we have
\begin{equation}\label{I-red}
\begin{aligned}
&\osc_{Q_{\frac1{8}\rho}(\theta)}u\le\bigg(1-\frac{\bar\eta}{|\ln(\bar\xi\om)|^{N-1}}\bigg)\om,\quad\text{where}\>\>\theta=\bigg(\frac{\om}{|\ln(\bar\xi\om)|^{N-1}}\bigg)^{2-N}.
\end{aligned}
\end{equation}

Now, we assume that the \textbf{second case} \eqref{II-alt} is satisfied. In the following two subsections, we perform the oscillation reduction for this case. Although the approach is very similar, we consider it in detail since now the regularization parameter $\eps$ plays a role.

To this end, we consider the following \textbf{two alternatives}:
\begin{align}
\exists\,\bar t\in I\quad\>\text{such that}\>\quad&|[u(\cdot,\bar t)\ge\mu^+-\tfrac18\om]\cap K_\rho|\le c_1(\tfrac18\om)^3|K_\rho|,\label{II-alt-i}\\
&\nonumber\\
\forall\,t\in I\quad\>\text{we have}\>\quad&|[u(\cdot,t)\ge\mu^+-\tfrac18\om]\cap K_\rho|> c_1(\tfrac18\om)^3|K_\rho|.\label{II-alt-ii}
\end{align}

\subsection{Second case -- first alternative}\label{SS:4:3}
Under the \textbf{second case} \eqref{II-alt},
we suppose the \textbf{first alternative} \eqref{II-alt-i} holds true, and we work precisely as we did in \S~\ref{SS:4:1} when we assumed that \eqref{I-alt-i} is satisfied. In particular, relying first on Lemma~\ref{Lm:DG-new} for sub-solutions, and then on Lemma~\ref{Lm:DG:initial:1} for sub-solutions, we obtain for any $\widetilde\xi\in(0,\tfrac1{32}]$  that 
\[
u\le\mu^+-\tfrac12\widetilde\xi\om
\]
in 
$$
\dsty K_{\frac1{8}\rho}
\times\Big(\bar t+\tfrac\delta{2^N}(\tfrac18\om)^{2-N}\rho^N,\bar t+\tfrac\delta{2^N}(\tfrac18\om)^{2-N}\rho^N+\gamma_o(\widetilde\xi\om)^{2-N}\tfrac1{4^N}\rho^N\Big],
$$ 
which corresponds to \eqref{Eq:lower-bd:1} and \eqref{Eq:cylinder:1}.
Afterwards, repeating almost verbatim the same computations, we choose $\widetilde\xi$ as in \eqref{choice-xi}, define $\theta$ as in \eqref{Eq:theta}, and conclude the same oscillation estimate as \eqref{osc-1}:
\[
\osc_{Q_{\frac1{8}\rho}(\theta)}u\le(1-\eta)\om,
\]
with $\theta$ as in \eqref{Eq:theta} and $\eta$ as in \eqref{theta-eta}.

\subsection{Second case -- second alternative}\label{SS:4:4}
Still under the \textbf{second case} \eqref{II-alt}, we consider instead the \textbf{second alternative} \eqref{II-alt-ii}. Since $\mu^-+\tfrac18\om\le\mu^+-\tfrac18\om$, we can rephrase \eqref{II-alt-ii} as
\[
\forall\,t\in I\quad\>\text{we have}\>\quad|[u(\cdot,t)\ge\mu^-+\tfrac18\om]\cap K_\rho|> c_1(\tfrac18\om)^3|K_\rho|.
\]
We introduce the new function
\[
v\df=\tfrac18\om-(u-k)_- , \quad\>\text{with}\>\,\,k=\mu^-+\tfrac18\om.
\]
Since $\mu^-+\tfrac18\om<-\tfrac18\om<0$, in order to apply Lemma~\ref{Lm:sub-solution}, and ensure that $v$ is a non-negative, weak super-solution to the parabolic equation~\eqref{Eq:3:1}, with structure conditions~\eqref{Eq:1:2} in $Q_o$, we need to require that $k<-\eps$; it suffices to assume that $-\tfrac18\om<-\eps$, that is, $\eps<\tfrac18\om$. Once such a condition on $\eps$ is satisfied, we can proceed almost verbatim as in \S~\ref{SS:4:2}, and conclude the same oscillation estimate as \eqref{osc-2}, \textit{i.e.},
\[
\osc_{Q_{\frac12\rho}(\theta)}u\le(1-\eta)\om,
\]
with $\theta$ as in\eqref{Eq:theta} and $\eta$ as in \eqref{theta-eta-bis}. 

Hence, under \eqref{osc-control} assumed at the beginning, we have eventually found that either $\om<8\eps$ or \eqref{I-red} holds true. Notice that \eqref{I-red} takes into account also the case when
$\mu^+-\mu^-<\tfrac12\om$. 

\subsection{Derivation of the modulus of continuity}
We have all the tools we need to derive a quantitative modulus of continuity.  
Let us summarize what we have achieved so far. Under the assumption that \eqref{osc-control} holds true for $\rho>0$, and $\om\in(0,1)$, we obtained that either
\begin{equation*}
\osc_{Q_{\frac18\rho}(\theta)}u\le\bigg(1-\frac{\bar\eta}{|\ln(\bar\xi\om)|^{N-1}}\bigg)\om,\quad\text{with}\>\,\,\theta=\bigg(\frac{\om}{|\ln(\bar\xi\om)|^{N-1}}\bigg)^{2-N},
\end{equation*}
where $\bar\eta$ depends only on the data $\{N,C_o,C_1\}$ and $\bar\xi$ depends additionally on $\nu$, 
or
\begin{equation*}
\om<8\eps.
\end{equation*}
Now, we need to iterate the argument. For that, we let
$$
\omega_o=\omega, \quad \rho_o=8\rho, \quad\theta_o=\theta,\quad\bar\theta_o=\bar\theta,
$$
set
\[
\begin{aligned}
& \omega_1\df=\bigg(1-\frac{\bar\eta}{|\ln(\bar\xi \omega_o)|^{N-1}}\bigg)\omega_o,\\
& \bar\theta_1\df=\Big[(\tfrac{1}{8} \omega_1) \bar{\delta} \operatorname{cap}_N\Big(K_{\bar c c_1^2(\frac{1}{8} \omega_1)^6}, K_3\Big)\Big]^{2-N},
\end{aligned}
\]
where $\bar\delta$, $\bar c$, $b$, $\xi$ are the quantities stipulated in Proposition~\ref{Prop:3:1} in terms of $\{N,C_o,C_1\}$, and $c_1=c_1(\nu,N,C_o,C_1)$ is the quantity in \eqref{omi},
and seek $\rho_1$ to verify the following set inclusion
\[
\begin{aligned}
& Q_{8\rho_1}( b \xi^{2-N} \bar\theta_1) \subseteq Q_{\frac{1}{64} \rho_o}(\theta_o).
\end{aligned}
\]
Without loss of generality, we may assume
\[
\frac{1}{\left|\ln(\bar\xi \omega_o)\right|^{N-1}} \le \frac{1}{2},
\]
which yields $\omega_1 \ge \frac{1}{2} \omega_o$. Hence, by the definition of $\bar\theta_1$, we can estimate
\[
8^N b\xi^{2-N} \bar\theta_1 \rho_1^N \le 8^N b \xi^{2-N}\Big[(\tfrac{1}{16} \omega_o) \bar{\delta}
\operatorname{cap}_N\Big(K_{\bar c c_1^2(\frac{1}{16} \omega_o)^6}, K_3\Big)\Big]^{2-N} \rho_1^N.
\]
Consequently, for the set inclusion, we select $\rho_1$ to satisfy
\[
\begin{aligned}
8^N b \xi^{2-N}\Big[(\tfrac{1}{16} \omega_o) \bar\delta
\operatorname{cap}_N\Big(K_{\bar c c_1^2(\frac{1}{16} \omega_o)^6}, K_3\Big)\Big]^{2-N} \rho_1^N
\le \omega_o^{2-N}|\ln(\bar\xi \omega_o)|^{(N-1)(N-2)}(\tfrac{1}{64}\rho_o)^N.
\end{aligned}
\]
To proceed, we rewrite the capacity term using \eqref{Eq:cap-dim:1}$_2$ and \eqref{Eq:cap-ln}:
\[
\begin{aligned}
\operatorname{cap}_N&\Big(K_{\bar c c_1^2(\frac{1}{16} \omega_o)^6}, K_3\Big)\\
& =b_2\bigg|\ln\Big[\frac{\bar c c_1^2}{2^6}\Big(\frac{\omega_o}{8}\Big)^6\Big]\bigg|^{-(N-1)}=b_2\bigg|\ln\Big(\frac{\bar\xi^6 \omega_o^6}{2^6}\Big)\bigg|^{-(N-1)} \\
& =b_2\big|\ln(\bar\xi^{6} \omega_o^6)-\ln 2^6\big|^{-(N-1)} =b_2\big[6\left|\ln(\bar\xi \omega_o)\right|+6\ln 2\big]^{-(N-1)} \\
& =b_2 6^{-(N-1)}\big[|\ln(\bar\xi \omega_o)|+\ln 2\big]^{-(N-1)}.
\end{aligned}
\]
Hence, plugging it into the last estimate yields
\[
\begin{aligned}
& 8^N b \xi^{2-N}  \big[(\tfrac{1}{16} \omega_o) \bar\delta 6^{-(N-1)} b_2\big]^{2-N}\big[|\ln(\bar\xi \omega_o) |+\ln 2\big]^{(N-1)(N-2)} \rho_1^N \\
&\qquad \le \om_o^{2-N} |\ln(\bar{\xi} \omega_o)|^{(N-1)(N-2)} \tfrac{1}{64^N} \rho_o^N.
\end{aligned}
\]
Without loss of generality, we may assume that
\[
|\ln(\bar\xi \omega_o)|>\ln 2\quad\Rightarrow\quad |\ln(\bar\xi \omega_o)|+\ln 2<2|\ln(\bar\xi 
\omega_o)|.
\]
Then, the above estimate is satisfied if
\[
8^N b \xi^{2-N} \big[\tfrac{1}{16} \bar\delta 6^{-(N-1)} b_2\big]^{2-N} \big[2|\ln(\bar\xi \omega_o)|\big]^{(N-1)(N-2)} \rho_1^N\le\tfrac{1}{64^N}|\ln(\bar\xi \omega_o)|^{(N-1)(N-2)} \rho_o^N.
\]
Hence, it suffices to choose $\rho_1$ such that
\[
\begin{aligned}
&8^N b \xi^{2-N}  2^{(N-1)(N-2)} \big[\tfrac{1}{16} \bar\delta 6^{-(N-1)} b_2\big]^{2-N}  \rho_1^N \le \tfrac{1}{64^N}  \rho_o^N.
\end{aligned}
\]
This suggests the choice 
\[
\rho_1=\lm\rho_o\qquad\text{with}\quad \lm\df=\frac{(\bar\dl \xi)^{\frac{N-2}{N}}}{\gm(N)b^{\frac1N}} .
\]
Note that $\lm$ depends only on the data $\{N,C_o,C_1\}$, and not on $\nu$. Hence, we can conclude that
\[
\osc_{Q_{8\rho_1}( b \xi^{2-N} \bar\theta_1)} u \le \omega_1,
\]
which at this stage plays the role of \eqref{osc-control}.

Repeating the arguments of Sections~\ref{SS:4:1}--\ref{SS:4:4} with the cylinder $Q_{8\rho_1}( b \xi^{2-N} \bar\theta_1)$, we obtain that either
\[
\osc_{Q_{\frac{1}{8} \rho_1}(\theta_1)} u \le\bigg(1-
\frac{\bar\eta}{|\ln(\bar{\xi} \omega_1)|^{N-1}}\bigg)\omega_1,\quad\text{ where }\,\,
\theta_1=\bigg(\frac{\om_1}{|\ln(\bar\xi\om_1)|^{N-1}}\bigg)^{2-N}
\]
or
\[
\omega_1 \le 8 \eps.
\]
Now, for every $n \in \nn_0$ we may construct 
\[
\begin{aligned}
& \rho_o=8\rho, \quad \rho_{n+1}=\lm\rho_n=\lm^{n+1}\rho_o=\lm^{n+1}8\rho, \\
& \omega_o=\omega, \quad \omega_{n+1}= \bigg(1-\frac{\bar\eta}{|\ln(\bar\xi \omega_n)|^{N-1}}\bigg) \omega_n  , \\
& \bar\theta_n=\bigg[(\tfrac{1}{8} \omega_n) \bar\delta \operatorname{cap}_N\Big(K_{\bar{c} c_1^2(\frac{1}{8} \omega_n)^6}, K_3\Big)\bigg]^{2-N}, \quad\theta_n=\bigg(\frac{\om_n}{|\ln(\bar\xi\om_n)|^{N-1}}\bigg)^{2-N},\\
& Q_n^{\prime}=Q_{8\rho_n}( b \xi^{2-N} \bar\theta_n), \quad Q_n=Q_{\frac{1}{8} \rho_n}(\theta_n).
\end{aligned}
\]

By induction, if up to some $j \in \nn$ we have
\[
\omega_n>8 \eps,\qquad\forall\, n \in\{0,1, \dots, j-1\},
\]
then, for all $n \in\{0,1,\dots, j\}$, there holds
\[
Q_n \subset Q_n^{\prime} \subset Q_{n-1}, \qquad \osc_{Q_n} u\le \osc_{Q_n^\prime} u \le \omega_n.
\]
On the other hand, we denote by $j$ the first index to satisfy
\begin{equation}\label{Eq:ast}
\omega_j \le 8 \eps.
\end{equation}
For $n\in\nn_0$, consider the sequence
\[
a_n=\exp (-c_* n^{\frac1N}), \quad\text{with}\quad c_* =\frac{\bar\eta}{2^{N-2}|\ln\bar\xi|^{N-1}}.
\]
In Appendix~\ref{A:2}, we prove that
\[
\left\{
\begin{aligned}
&a_{n+1} \ge a_n\bigg(1-\frac{\bar\eta}{|\ln(\bar\xi a_n)|^{N-1}}\bigg), \\
&a_o \ge \omega_o.
\end{aligned}
\right.
\]
Hence, for any $ n \in\nn_0$, we have $\displaystyle a_n \ge \omega_n$. 

Let us now take $r \in(0, \rho)$; if, for some $n \in\{0,1,\dots,j\}$, we have
$\displaystyle \rho_{n+1} \le  8r < \rho_n$, then 
\[
\osc_{Q_r(\theta_o)} u\le \osc_{Q_n} u \le\omega_n
\le \exp \big(-c_* n^{\frac1N}\big).
\]
On the other hand, we have $\displaystyle \rho_{n+1}=\lm^{n+1}\rho$. 
Therefore, $\rho_{n+1} \le  8r$ yields
\[
\begin{aligned}
 n \ge \frac{1}{2|\ln \lm|} \left| \ln\left(\frac{8r}{\rho}\right) \right|,
\end{aligned}
\]
and we conclude that
\[
\begin{aligned}
\osc_{Q_r(\theta_o)} u&\le
\exp\left(-\frac{c_*}{2^{\frac1N}|\ln\lm|^{\frac1N}}\left|\ln\left(\frac{8r}{\rho}\right)\right|^{\frac1N}\right)\\
&=\exp\left(-c\left|\ln\left(\frac{8r}{\rho}\right)\right|^{\frac1N}\right),
\quad\mbox{where $\displaystyle c\df=\frac{c_*}{2^{\frac1N}|\ln\lm|^{\frac1N}}$.}
\end{aligned}
\]
It is not hard to trace the dependence of $c$ on $\nu_* = \max\{1,\nu\}$ as
\begin{equation}\label{Eq:c}
    c=\frac{\gm(N,C_o,C_1)}{\ln (2\nu_*)^{N-1}}.
\end{equation}
Therefore, if we know $\nu_*$ is bounded by a larger number $ \bar\nu$, then we may replace $\nu_*$ by $\bar\nu$ in the definition of $c$ while retaining the same oscillation estimate.

On the other hand, if $ r<\rho_{j+1}$, where $j$ is the first index for which \eqref{Eq:ast} holds, we may use it and conclude that
\[
\osc_{Q_r(\theta_o)} u\le\osc_{Q_j} u\le\omega_j\le8\eps.
\]
This allows us to include the $\eps$-term into the oscillation estimate and obtain that
\begin{equation*}
\begin{aligned}
\osc_{Q_r(\theta_o)} u \le\exp\left(-c\left| \ln\left(\frac{8r}{\rho}\right)\right|^{\frac1N}\right)+8 \eps.
\end{aligned}
\end{equation*}
Finally, we can let $\eps \rightarrow 0$ and conclude with the wanted modulus of continuity.

\subsection{Modulus of continuity over compact sets} \label{4.6}
The gist of the previous arguments is that if there are parameters $\rho>0$ and $\om_o\in(0,1)$ such that the initial oscillation estimate \eqref{osc-control} is verified, then the oscillation decay
\[
\osc_{Q_r(\theta_o)}u\le \exp\left(-c\left|\ln\left(\frac{8r}{\rho}\right)\right|^{\frac1N}\right)
\]
 can be reached for all $r\in(0, \rho)$.
 The initial oscillation estimate \eqref{osc-control} can be rewritten as
 \begin{equation}\label{osc-control-1}
     \osc_{Q_{8\rho}(L\theta_o)} u\le \om_o,
 \end{equation}
 where $L>1$ can be easily calculated in terms of the data $\{N, C_o, C_1\}$, making use of \eqref{Eq:cap-ln}.
Moreover, the condition $\om_o\le1$ yields
\[
1<\theta_o=\left(\frac{\om_o}{|\ln(\bar\xi \om_o)|^{N-1}}\right)^{2-N}, 
\]
and hence, under the assumption \eqref{osc-control-1}, we conclude that 
\begin{equation}\label{Eq:osc:final}
\osc_{Q_r}u\le \exp\left(-c\left|\ln\left(\frac{8r}{\rho}\right)\right|^{\frac1N}\right).
\end{equation}
We want to show how such a modulus of continuity is affected as we approach the parabolic boundary of $E_T$.

Let $\mathcal K$ be a compact subset of $E_T$, and define 
\[
\mathbb{M}\df=\max\Big\{1,\,\osc_{E_T}u\Big\},\qquad 8R\df=p-\dist_{\textrm{par}}({\mathcal K},\partial_{\textrm{par}}E_T).
\]
By the definition of $p-$parabolic distance, $\forall\, (x_o,t_o)\in{\mathcal K}$, we have
\[
(x_o,t_o)+Q_{8R}(\mathbb{M}^{2-N})\subset E_T. 
\]
Fix such a point $(x_o,t_o)\in{\mathcal K}$ and consider a new function on $Q_{8R}$ defined by
\[
\displaystyle v(x, t)\df=\frac {u(x-x_o, \mathbb{M}^{2-N}(t-t_o))}{\mathbb M}.
\]
Then, apparently, we have
\[
\osc_{Q_{8R}}v\le 1.
\]
Moreover, such $v$ satisfies the same kind of Stefan problem \eqref{Eq:1:1}, with structure conditions \eqref{Eq:1:2} but with a different jump constant $\nu/\mathbb M$ that defines $\be(\cdot)$ in \eqref{Eq:beta}.
The last oscillation estimate of $v$ yields an analog of \eqref{osc-control-1} if we let $\om_o=1$ and $\rho=\sig R$, where 
\[
\sig=\frac1{L^{\frac1N}|\ln\bar\xi|^{\frac{(N-2)(N-1)}{N}}}
\]
depends on $\{\nu,N,C_o,C_1\}$.
Therefore, we may reproduce all previous arguments for $v$ now.
Moreover, due to $\mathbb M\ge1$, the jump constant $\nu/\mathbb M$ that appears in all estimates regarding $v$ can be replaced by the larger constant $\nu$, cf.~\eqref{Eq:c}.
As a result, we obtain an analog of oscillation decay \eqref{Eq:osc:final} for $v$, with the constant $c$ still determined by the data $\{\nu, N, C_o, C_1\}$. 
Reverting to $u$, we arrive at
\[
\osc_{(x_o,t_o)+Q_r(\mathbb M^{2-N})} u\le \mathbb M \exp\left(-c\left|\ln\left(\frac{8r}{\sig R}\right)\right|^{\frac1N}\right),
\]
for any $r\in (0, \sig R)$. 

Next, pick a second point $(x_1,t_1)\in{\mathcal K}$. Assume with no loss of generality that $t_1\le t_o$ and
\[
r\df= |x_o-x_1|+\mathbb M^{\frac{N-2}{N}}(t_o-t_1)^{\frac1N}< \sig R,
\]
such that
\[
(x_1,t_1)\in (x_o,t_o)+Q_{r}(\mathbb M^{2-N}).
\]
The oscillation decay then yields
\[
|u(x_o,t_o)-u(x_1,t_1)|\le \mathbb M \exp\left(-c\left|\ln\left(\frac{|x_o-x_1|+\mathbb M^{\frac{N-2}{N}}(t_o-t_1)^{\frac1N}}{\frac18\sig R} \right)\right|^{\frac1N}\right).
\]
This is the desired modulus of continuity over the compact set $\mathcal{K}$.

\subsection{H\"older modulus of continuity}
In this section, we briefly indicate how to modify the argument of the last subsections in order to obtain the H\"older continuity of weak solutions in the case $p>\max\{2,N\}$.

As before, the starting point is an oscillation estimate like \eqref{osc-control}.
Indeed, suppose there are parameters $\rho>0$ and $ \om\in(0,1)$ such that
\begin{equation}\label{osc-control-2}
Q_{8\rho}(b\xi^{2-p}\bar\theta)\subset E_T,\qquad\osc_{Q_{8\rho}(b\xi^{2-p}\bar\theta)} u\le\om,
\end{equation}
where
\begin{equation}\label{al-theta-2}
\bar\theta\df=\left[\left(\tfrac18\om\right)\bar\delta\,b_1\right]^{2-p},\qquad\al\df=c_1\left(\tfrac18\om\right)^{\frac{N+2p}{p}}.
\end{equation}
The various parameters, $b$, $c_1$, $\bar\dl$, $\xi$, etc... retain the same token as in \eqref{osc-control} and \eqref{al-theta}. However, for $\bar\theta$ in \eqref{al-theta-2}, in the place of the capacity, we employ its lower bound $b_1$ from \eqref{Eq:cap-bd}, which depends only on $\{N,p\}$. For $\al$, we change the power in accordance with \eqref{omi}.

Next, we introduce $\mu^{\pm}$ and consider, analogously, the two cases \eqref{I-alt} and \eqref{II-alt}. Clearly, we define the interval $I$ with the parameter $N$ in the power replaced by $p$, whereas the two alternatives \eqref{I-alt-i} and \eqref{I-alt-ii} remain the same upon taking these changes into consideration. 

Section~\ref{SS:4:1} can be reproduced along the same lines. Clearly, the parameter $N$ that appears in the power of various quantities must be changed to $p$. As a result, an analog of the pointwise estimate \eqref{Eq:lower-bd:1} is reached in the cylinder \eqref{Eq:cylinder:1}. This will yield an oscillation estimate provided we choose $\widetilde{\xi}$ and $\xi$ properly, to guarantee 
\begin{align}\label{Eq:t-bar:1-}
    \tfrac\delta{2^p}(\tfrac18\om)^{2-p}\rho^p+\gamma_o(\widetilde\xi\om)^{2-p}\tfrac1{4^p}\rho^p&\ge 2b\xi^{2-p}\bar\theta \rho^p, \\\label{Eq:t-bar:2-}
    -b\xi^{2-p}\bar\theta\rho^p+\tfrac\delta{2^p}(\tfrac18\om)^{2-p}\rho^p&<-\theta(\tfrac18\rho)^p,
\end{align}
where 
\begin{equation}\label{Eq:theta-om}
    \theta=\om^{2-p}.
\end{equation}
Discarding the term with $\dl$ in \eqref{Eq:t-bar:1-} and considering $\bar\theta$ defined in \eqref{al-theta-2}, the choice of $\widetilde{\xi}$ now becomes
\[
\widetilde\xi\le \Big(\frac{\gamma_o}{b}\Big)^{\frac1{p-2}}2^{\frac{5(1-p)}{p-2}} \xi\, \bar\delta\, b_1,
\]
whereas \eqref{Eq:t-bar:2-} is satisfied provided $\xi$ is small enough. 
Therefore, we arrive at an analog of the oscillation estimate \eqref{osc-1}, with $\theta$ defined in \eqref{Eq:theta-om} and with $\eta$ replaced by $\frac12\widetilde{\xi}$, which depends only on the data $\{N,p,C_o,C_1\}$ and is independent of $\om$ and $\nu$.

Section~\ref{SS:4:2} can also be reproduced along the same lines. However, in the application of Proposition~\ref{Prop:3:1} to $v$, we obtain instead that
\[
v\ge\tfrac18\om\,\xi\,b_1
\]
in $K_{\frac12\rho}\times(-b\xi^{2-p}\bar\theta\rho^p,0]$, because of the capacity estimate in \eqref{Eq:cap-bd}. Consequently, we obtain an analog of the oscillation estimate \eqref{osc-2}, with $\theta$ defined in \eqref{Eq:theta-om} and with $\eta$ replaced by $\frac18\xi b_1$, which depends only on the data $\{N,p,C_o,C_1\}$, and is independent of $\om$ and $\nu$.

Necessary changes can be performed similarly for Sections \ref{SS:4:3} and \ref{SS:4:4}. Therefore, under the assumption that \eqref{osc-control-2} holds true for $\rho>0$ and $\om\in(0,1)$, we eventually obtain
\begin{equation*}
    \osc_{Q_{\frac18\rho}(\theta)} u\le (1-\eta)\om,\quad\text{with}\>\>\theta=\om^{2-p}.
\end{equation*}
Note that $\eta=\min\{\frac12\widetilde{\xi},\frac18 \xi b_1\}$ depends only on $\{N,p,C_o,C_1\}$, but neither on $\om$ nor on $\nu$.
Setting up an iteration scheme and deriving an H\"older modulus of continuity now become standard. 

Finally, it is apparent that when $2<p<N$, $N\ge3$, Proposition~\ref{Prop:3:1} yields a power-like dependence on $\al$ of the pointwise estimate. This will affect the recurrence of $\omega_n$ and, consequently, the modulus of continuity. Once this is considered, the previous computations can also be appropriately adapted for this case. All the moduli of Table~\ref{Table:1} are now justified.

\appendix

\section{Some properties of $p-$capacity}\label{A:1}

For $p\ge1$, the notion of $p-$capacity is defined by
\[
{\rm cap}_p(F, \Om)\df{=}\inf_{u\in W}\int_\Om |Du|^p\,\dx,
\]
where $F$ is a compact subset of the open set $\Om$ in $\rn$
and 
$$
W\df=\{u\in C_o^{\infty}(\Om):u\ge1\>\text{on}\>F\}.
$$

After proper approximations, the above minimization could take place over $W^{1,p}_o(\Om)$ instead of $C_o^{\infty}(\Om)$,
and the $p-$capacity of an arbitrary subset of $\Om$
can be formulated based on that of compact subsets. The reader is referred to \cite[Chapter~2]{HKM}.

Immediate from the definition is the scaling property
\begin{equation}\label{Eq:cap-scale}
\frac{{\rm cap}_p(K_{\varep\rho}, K_{\rho})}{\rho^{N-p}}={\rm cap}_p(K_{\varep}, K_{1}),\quad\text{for}\>\varep\in(0,1).
\end{equation}
Also immediate is that 
$$
{\rm cap}_p(F, \Om)\ge{\rm cap}_p(F, \Om_o), 
$$ 
if $\Om_o$ is an open set satisfying $F\subset\Om\subset\Om_o$. This, in particular, implies that 
$$
{\rm cap}_p(K_{\varep}, K_{1})\ge{\rm cap}_p(K_{\varep}, K_{3}).
$$ 
With a little more effort, one also shows that the reverse estimate holds, apart from a multiplicative constant depending on $N$; see \cite[\S~2.16]{HKM}.

Capacities of balls can be estimated explicitly. Namely, 
\begin{equation}\label{Eq:cap-dim:0}
{\rm cap}_p(B_{r}, B_{R}) = \left\{
\begin{array}{cl}
\om_{N-1}\left(\frac{|N-p|}{p-1}\right)^{p-1} \left|R^{\frac{p-N}{p-1}} - r^{\frac{p-N}{p-1}}\right|^{1-p} &\quad \text{if}\>p\neq N,\\[5pt]
\om_{N-1}\left|\ln\frac{r}{R}\right|^{1-N} &\quad \text{if}\>p=N, 
\end{array}\right.
\end{equation}
where $\om_{N-1}$ denotes the surface measure of the boundary of the unit ball in $\rn$. See \cite[\S~2.11]{HKM}. 
Based on this, 
one easily estimates
\begin{equation}\label{Eq:cap-bd}
{\rm cap}_p(K_{\varep}, K_{3})\le \gm(N,p),
\end{equation}
for any $\varep\in(0,1)$. 

Another fact that can be derived from \eqref{Eq:cap-dim:0} with little effort is 
\begin{equation}\label{Eq:cap-dim:1}
{\rm cap}_p(K_{\varep}, K_{3}) \,\ge \left\{
\begin{array}{lr}
b_1 &\quad \text{if}\>p>N,\\[5pt]
b_2|\ln\varep|^{-(N-1)} &\quad \text{if}\>p=N,\\[5pt]
b_3 \varep^{N-p} &\quad \text{if}\>p<N,
\end{array}\right.
\end{equation}
for positive constants $b_1$, $b_2$ and $b_3$ depending only on $\{N,\,p\}$, and the lower bound is actually an equality when $p=N$. 

\section{Estimate of a parameter}\label{A:2}
Consider the decreasing sequence 
\[
a_n=\exp(-c_*\, n^{\frac1N}),\qquad n\in\nn_0.
\]
Our goal is to select the positive parameter $c_*$ in such a way that, for any $n\in\N_0$,
\begin{equation}\label{def-an}
\left\lbrace
\begin{aligned}
&a_{n+1}\ge a_n\bigg(1-\frac{\bar\eta}{|\ln(\bar\xi a_n)|^{N-1}}\bigg),\\
&a_o\ge\om,
\end{aligned}
\right.
\end{equation}
where $\bar\eta$, $\bar\xi$, and $\om$ are the same quantities as in \eqref{I-red}.

Since $a_o=1$, and in the derivation of \eqref{I-red}, we assumed $\om\le1$, the second condition of \eqref{def-an} is automatically satisfied.

Moreover, we have
\begin{align*}
1-\frac{\bar\eta}{|\ln(\bar\xi a_n)|^{N-1}}&=1-\frac{\bar\eta}{[|\ln a_n|+|\ln \bar\xi|]^{N-1}}\\
&=1-\frac{\bar\eta}{[c_*\,n^{\frac1N}+|\ln \bar\xi|]^{N-1}}
\le\exp\bigg(-\frac{\bar\eta}{[c_*\,n^{\frac1N}+|\ln \bar\xi|]^{N-1}}\bigg),
\end{align*}
where we have taken into account that one has $1-x\le e^{-x}$, $\forall\,x\in\rr$.

On the other hand, since
\[
\frac{a_{n+1}}{a_n}=\exp\Big(-c_*[(n+1)^{\frac1N}-n^{\frac1N}]\Big),
\]
it suffices to choose $c_*$ such that 
\[
c_*[(n+1)^{\frac1N}-n^{\frac1N}]\le\frac{\bar\eta}{[c_*\,n^{\frac1N}+|\ln \bar\xi|]^{N-1}},
\]
that is
\[
c_*\big[(n+1)^{\frac1N}-n^{\frac1N}\big]\big[c_*\,n^{\frac1N}+|\ln \bar\xi|\big]^{N-1}\le \bar\eta.
\]
Since, for any $a,\,b>0$, we have $(a+b)^{N-1}\le 2^{N-2}(a^{N-1}+b^{N-1})$, we require
\[
\begin{aligned}
2^{N-2} c_* \big[(n+1)^{\frac1N}-n^{\frac1N}\big]\big[c_*^{N-1}\,n^{\frac{N-1}N}+|\ln \bar\xi|^{N-1}\big]&\le \bar\eta\\
c_*^N\big[(n^N+n^{N-1})^{\frac1N}-n\big]+c_*\big[(n+1)^{\frac1N}-n^{\frac1N}\big]|\ln\bar\xi|^{N-1}&\le\frac{\bar\eta}{2^{N-2}}.
\end{aligned}
\]
If $n=0$, we require 
$$
\displaystyle c_*\le\frac{\bar\eta}{2^{N-2}|\ln\bar\xi|^{N-1}}.
$$ 
If $n\ge1$, recalling that, for any $m\in(0,1)$ and for any $x>-1$, we have $(1+x)^m\le 1+mx$, yields
\[
(n^N+n^{N-1})^{\frac1N}-n=n\bigg(1+\frac1n\bigg)^{\frac1N}-n\le n\bigg(1+\frac1N\frac1n\bigg)-n=\frac1N,
\]
\[
(n+1)^{\frac1N}-n^{\frac1N}=n^{\frac1N}\bigg(1+\frac1n\bigg)^{\frac1N}-n^{\frac1N}\le n^{\frac1N}+\frac1N\frac1{n^{\frac{N-1}N}}-n^{\frac1N}\le\frac1N.
\]
Hence, we require
\begin{align*}
c_*^N+c_*|\ln\bar\xi|^{N-1}&\le\frac{N \bar\eta}{2^{N-2}},\\
2c_*|\ln\bar\xi|^{N-1}&\le\frac{N \bar\eta}{2^{N-2}},\\
c_*&\le\frac{N \bar\eta}{2^{N-1}|\ln\bar\xi|^{N-1}}.
\end{align*}
Therefore, in order to have $c_*$ that satisfies \eqref{def-an} for any $n\ge0$, it suffices to take
\[
c_*\df=\min\left\{\frac{N \bar\eta}{2^{N-1}|\ln\bar\xi|^{N-1}},\frac{\bar\eta}{2^{N-2}|\ln\bar\xi|^{N-1}}\right\}=\frac{\bar\eta}{2^{N-2}|\ln\bar\xi|^{N-1}}.
\]

\bigskip

{\small \noindent{\bf Acknowledgments.} UG is partially supported by the grant 202244A7YL ``Gradient flows and non-smooth geometric structures with applications to optimization and machine learning." NL is supported by the FWF-project P36272-N ``On the Stefan type problems." JMU is partially supported by KAUST and CMUC (funded by the Portuguese Government through FCT/MCTES, DOI 10.54499/UIDB/00324/2020). Part of this paper was written while the authors participated in the Workshop on Nonlinear Parabolic PDEs, held at the Institut Mittag-Leffler in May 2024.}

\end{document}